\documentclass[11pt,reqno]{amsart}
\usepackage{graphicx} % Required for inserting images
\usepackage{comment}
\usepackage[pagebackref, colorlinks=true,linkcolor=red,citecolor=blue]{hyperref}
\usepackage{latexsym}
\usepackage{amsmath}
\usepackage{amsfonts}
\usepackage{amssymb}
\usepackage{tikz}
\usepackage{amsthm}
\usetikzlibrary{arrows, shapes}
\usepackage{graphicx}
\usepackage{enumitem}
\usepackage{mathrsfs}
\usepackage{enumitem}
\usepackage{hyperref}
\usepackage[arrow, matrix, curve]{xy}
\usepackage[nameinlink,capitalize]{cleveref}
\usepackage{caption}
\usepackage{subcaption}
\usepackage{color}


 
\newtheorem{theorem}{Theorem}[section]
\newtheorem{lemma}[theorem]{Lemma}
\newtheorem{corollary}[theorem]{Corollary}
\newtheorem{proposition}[theorem]{Proposition}

\theoremstyle{definition}
\newtheorem{definition}[theorem]{Definition}
\newtheorem{remark} [theorem]{Remark}
\newtheorem{example}[theorem]{Example}

\newcounter{manualthm}

\theoremstyle{plain}
\newtheorem{manualtheorem}[manualthm]{Theorem}

\newenvironment{manualtheoremx}[1]{%
  \refstepcounter{manualthm}% step the counter so \label works
  \renewcommand{\themanualthm}{#1}% override the number
  \begin{manualtheorem}%
}{\end{manualtheorem}}

\theoremstyle{plain}

\theoremstyle{remark}
\newtheorem{Remark}[theorem]{\rm \bf Remark}
\numberwithin{equation}{section}

\newcommand{\hg}[1]{{}_2F_1\Big( \frac{#1}{2}\Big)}
\newcommand{\hgf}[1]{{}_2F_1( #1 )}

\newcommand{\ltfrac}[2]{\mbox{\Large$\frac{#1}{#2}$}}
\newcommand{\mr}[1]{\mathrm{#1}}

\newcommand{\cl}[1]{\mathrm{cl}\, (#1)}

\newcommand{\C}[1]{\mathbb{C}^{#1}}
\newcommand{\bbC}{\mathbb{C}}
\newcommand{\D}{\mathbb{D}}
\newcommand{\R}[1]{\mathbb{R}^{#1}}
\newcommand{\bbR}{\mathbb{R}}
\newcommand{\N}{\mathbb{N}}

\newcommand{\bx}{\mathbf{x}}
\newcommand{\bv}{\mathbf{v}}

\newcommand{\by}{\mathbf{y}}

\newcommand{\fa}{\mathfrak{a}}
\newcommand{\fg}{\mathfrak{g}}
\newcommand{\fh}{\mathfrak{h}}
\newcommand{\fk}{\mathfrak{k}}

\newcommand{\fm}{\mathfrak{m}}
\newcommand{\fn}{\mathfrak{n}}

\newcommand{\fp}{\mathfrak{p}}

\newcommand{\fq}{\mathfrak{q}}
\newcommand{\fz}{\mathfrak{z}}

\newcommand{\cD}{\mathcal{D}}
\newcommand{\cT}{\mathcal{T}}
\newcommand{\rC}{\mathrm{C}}

\newcommand{\rJ}{\mathrm{J}}
\newcommand{\rI}{\mathrm{I}}
\newcommand{\rO}{\mathrm{O}}
\newcommand{\rS}{\mathrm{S}}
\newcommand{\rT}{\mathrm{T}}
\newcommand{\Supp}{\mathrm{Supp\, }}
\newcommand{\oz}{\overline{z}}
\newcommand{\ow}{\overline{w}}

 \newcommand{\oPsi}{\overline{\Psi}}
\newcommand{\SO}{\mathrm{SO}}
\newcommand{\GL}{\mathrm{GL}}

\newcommand{\so}{\mathfrak{so}}
\newcommand{\ad}{\mathrm{ad\, }}
\newcommand{\Hn}{\mathbb{H}^n}
\newcommand{\oHn}{\overline{{\mathbb H}}^n}
\newcommand{\dS}{\mathrm{dS}^n}
\newcommand{\dSC}{\mathrm{dS}^n_{\mathbb{C}}}

\newcommand{\Sn}{\mathbb{S}^n}
\newcommand{\oXi}{\overline{\Xi}}
\newcommand{\Scut}{(\dS \times \dS)_{\text{cut}}}
\newcommand{\Exp}{\mathrm{Exp}}
\newcommand{\Fl}{{}_2F_1}

\let\oldproofname=\proofname
\renewcommand{\proofname}{\rm\bf{\oldproofname}}
\newcommand{\ip}[2]{\langle #1,#2\rangle} 
\renewcommand{\:}{\, : \,}

\newcommand{\sgn}[1]{\mathrm{sgn}( #1)}

\newcommand{\DdS}{\mathcal{D}^*(\dS)}
\newcommand{\DdSl}{\mathcal{D}^*_\lambda(\dS)}

\title[Analytic wavefront set on De Sitter space]{Analytic wavefront sets of spherical distributions on de Sitter space}

\author{Gestur \'{O}lafsson}
\address{Department of Mathematics, Louisiana State University, Baton Rouge, LA 70803, USA}
\email{olafsson@math.lsu.edu}

\author{Iswarya Sitiraju}
\address{Department of Mathematics, Indian Institute of Technology, Bombay, MH 400076, India}
\email{iswarya@math.iitb.ac.in}

\begin{document}
\begin{abstract}
In this article we determine the wavefront sets of  spherical distributions on the de Sitter space 
$\dS =G/H$, $G= \SO_{1,n}(\R{})_e$. These are eigendistributions of the Laplacian on $\dS = G/H$ invariant
under the subgroup $H$. We construct bases for the spaces of  
spherical distributions as boundary values of sesquiholomorphic kernels on a certain $G$-invariant complex domain in
$\dSC$   containing the de Sitter space  as a $G$-orbit on the boundary. We characterize the elements of the basis by their 
 analytic wavefront sets. We also treat the spherical distributions invariant under $\mr{O}_{1,n-1}(\R{})$.
\end{abstract}

\maketitle

\tableofcontents 

\section{Introduction} Let $X$ be a manifold and $\Theta $ a distribution on $X$, i.e., a continuous
conjugate linear functional $C^\infty_c(X) \to \C{}$. 
The {\it wavefront set} $WF(\Theta )$ of $\Theta$ is a subset of $i\rT^*(X)$ encoding information about the singularities of the distribution. The
projection onto $X$ is the singular support, while the projection onto $i\rT_x(X)$, $x\in X$, gives the directions in which the
singularities at $x$ occur.  If the manifold $X$ is analytic, then the singularities of  $\Theta$ can also be 
studied in terms of the {\it analytic wavefront set}  $WF_A(\Theta)$.  This is the set containing the smooth wavefront set $WF(\Theta)$ and
consisting of points having no neighborhood on which 
$\Theta$ is given by a real analytic function.   

The term wavefront set was introduced in 1963 by Lars H\"ormander \cite[Chap. 8]{H63}. Later on it was used by the same author in \cite{H70,H71}; see also \cite{H03}. One of the  ideas was 
to use the microlocal information contained in the wavefront set to study the propagation of singularities of 
pseudo-differential operators. The wavefront set has evolved into a crucial concept in microlocal analysis and
the study of pseudo-differential operators. As such, it is an indispensable tool  in quantum field theory (QFT). 
One of the earliest papers to use wavefront sets in QFT is \cite{Di79}; later works include \cite{RM,SVW02, V99}. 

The first to extend the idea to representation theory was  R. Howe in \cite{HR}. We refer to \cite{BW25,HHO} for
 more recent results and references\ in this direction and to \cite{HO24} for an application to the
 Plancherel Theorem of spherical homogeneous spaces.
 
 In this article, we  focus on the wavefront set and the  analytic wavefront set of
 spherical  distributions  on the de Sitter space. They play a central role in 
 harmonic analysis on the de Sitter space as well as in quantum field theory on this space. The two main results are (1) 
 an explicit description of the wavefront set of a natural basis for the
 space of spherical distributions and (2) a proof that the analytic  wavefront set characterizes those distributions. 

We write elements in $\R{n+1}$ as $x=(x_0,x_1,\ldots ,x_n) = (x_0,\bx)$. For $n \geq 2$, let $[\cdot,\cdot]$ be the Lorentzian bilinear form  
\[[x,y] = -x_0y_0 +x_1y_1+\cdots +x_ny_n = -x_0y_0 +\ip{\bx}{\by}.\]
We denote by $\R{1,n}$ the space $\R{n+1}$ with the Lorentzian form $[\, \cdot\,  ,\, \cdot\,  ]$. It is
the $(n+1)$-dimensional Minkowski space.
 The de Sitter space is the one-sheeted hyperboloid defined by
\[  \dS= \{x \in \R{1,n}\: [x,x] = 1 \}.\]
It is one of the simple models of the universe in general relativity.
The de Sitter space is contained in the complexification 
\[ \dSC = \{z \in \C{n+1}\: [z,z]=1\}.\]

For our purposes, the de Sitter space can be realized as a homogeneous space in two different ways. Let 
$G^\prime$ be the isometry group of the Lorentzian form:
\[G' = \rO_{1,n}(\R{})= \{ g\in \GL_{n+1}(\R{})\: (\forall x,y\in \R{1,n})\, [gx,gy] = [x,y]\}\]
and let $G = \mr{SO}_{1,n} (\R{})_e$ be the connected component of $G^\prime$ containing the identity element.  
Let $H' \subset G'$ and $H \subset G$ be the stabilizers of $e_n = (0,...,0,1)$. Then
$H^\prime \simeq \rO_{1,n-1}(\R{})$ and $H \simeq \SO_{1,n-1}(\R{})_e$ are closed subgroups and
\[\dS \simeq G'/H' \simeq G/H. \]

Denote by $\Omega$ the future light cone, 
\[\Omega := \{v \in \R{1,n}\: [v,v] < 0, v_0 > 0\}\simeq \mathbb{R}^+\SO_{1,n-1}(\R{})_e\cdot e_0\]
and the corresponding
 tube domain by
\[\rT_\Omega := \{u + iv\: u,v \in \R{1,n}, v \in \Omega\}.\] It was shown in \cite{NO20} that the de Sitter space lies on the 
boundary of the  complex domain $\Xi=\dSC\cap \rT_\Omega$.
It was also shown that $\Xi$   is
the crown domain associated to the hyperboloid 
\[\Hn=\{ix\in i\R{1,n}\: -[x,x]=1\}=G\cdot (ie_0) \subset i\R{1,n}\cap \dSC, \]
 where $ e_0 =(1,0,\ldots ,0) $; see  \cite{AG90,KSt04} and the references therein for the definition and discussion
 of the crown. Furthermore, the authors of \cite{NO20} considered the   set of positive-definite $G$-invariant sesquiholomorphic kernels 
 $\Psi$ on $\Xi \times\Xi$ normalized by $\Psi(ie_0,ie_0)=1$ and showed that
 they are given by hypergeometric functions and parametrized by $i[0,\infty)\cup (0,(n-1)/2)$; see below for
 more details.

In this work, we study the kernels $\Psi_{\lambda}$, $\lambda \in \C{}$, introduced in
\cite{NO20} and normalized by $\Psi_\lambda (ie_0,ie_0)=1$. It was shown in \cite[Thm. 4.12]{NO20} and \cite[Cor. 4.13]{NO20}
 that they are given by
\begin{equation}\label{eq:1.1}  
  \Psi_{\lambda}(z,w) =  {}_2F_1\Big(\rho+\lambda,\rho-\lambda; \frac{n}{2}; \ltfrac{1+[z,\Bar{w}]}{2}\Big), \quad  z,w \in \Xi, 
\end{equation}
where we have simplified the notation by setting $\rho = (n-1)/2$. They are $G$-invariant sesquiholomorphic kernels on
$\Xi\times \Xi$ satisfying $\Psi_\lambda = \Psi_{-\lambda}$. Furthermore, the kernel $\Psi_\lambda$ is positive
definite if and only if $\lambda \in i\R{} \cup (-\rho ,\rho)$ (see \cite[Cor. 5.11]{NO20}). Moreover, for a fixed $w\in\Xi$, the
function $z\mapsto \Psi_\lambda (z,w)$ is an eigenfunction of the Laplacian $\Delta$:
\begin{equation}\label{eq:1}
\left(\Delta - \left(\rho ^2-\lambda^2\right)\right)\Psi_\lambda (\cdot ,w) = 0,
\end{equation}
and is invariant under the stabilizer group of $w$. Finally, the function $x\mapsto \Psi_\lambda (x,ie_0)$, $x\in \Hn$, is the
spherical function on $\Hn$ with spectral parameter $\lambda$; see \cite[Thm. 5.10]{NO20}. A celebrated theorem of Helgason and Johnson \cite{HJ69} implies that the spherical function with spectral parameter $\lambda$ is bounded on $\Hn$ if and only if $\lambda$ belongs to the strip \[[-\rho,\rho] \oplus i\R{} = \{z\in \C{}\: \Re z \in [-\rho,\rho]\}.\] But it is positive definite only if $\lambda \in \R{}$ or $ \lambda \in (-\rho ,\rho)$; see \cite[Thm. 7.5.9]{D09}.

 %We also need the conjugate kernels on $\oXi$ given by the same formula, but defined on a different domain:
%\begin{equation}\label{eq:pm2}
 %   \ti{\Psi}_{\lambda}({z},{w})   =
  %   {}_2F_1\Big(\rho+\lambda,\rho-\lambda; \frac{n}{2}; \ltfrac{1+[{z},\ow]}{2}\Big) \quad  {z},{w} \in \overline{\Xi}
%\end{equation}
% Note that $\oPsi_\lambda (z,w) = \Psi_\lambda (\oz,\ow)$ for
%$z,w\in\oXi$.

For each $z\in \Xi$, the function 
\begin{equation}\label{eq:1.2}
y\mapsto\Psi_{\lambda ,z}(y) = \Psi_\lambda (z,y)
\end{equation}
is analytic and we show that the boundary value   $\Xi\ni z\to x\in \dS$  exists in the space
of distributions on $\dS$; see Theorem 1 and the discussion in  Sections 2 and 3 for details.  We denote
the limiting distribution by  $\Psi_{\lambda,x} $. By \eqref{eq:1}, this is an eigendistribution of the Laplace-Beltrami operator $\Delta$ on 
$\dS$,  and $\Psi_{\lambda ,e_n}$ is $H$-invariant and hence a {\it spherical distribution}.
Moreover, the distribution $\Psi_{\lambda,x}$ is real analytic on the cut domain
\begin{equation}\label{eq:1.3}
\Scut = \dS \times \dS \setminus \{(x,y): [x-y,x-y] =0 \}
\end{equation}
and has 
jump discontinuities along the cut. These singularities are studied in terms of the analytic wavefront sets, which
we describe explicitly. We obtain from Theorem~\ref{thm:supp} that the support of the distribution $\Psi_{\lambda,x}$  is the entire space $\dS$.

The complexification $\dSC$ is invariant under complex conjugation $z\mapsto \oz$ and hence
$\oXi =\{\oz \mid z\in \Xi\}\subset \dSC$. In fact 
 \[\oXi =\dSC\cap \overline{\rT_\Omega} = \dSC\cap \rT_{-\Omega},\]
 which is the crown of the Riemannian hyperbolic space realized as $\oHn:=G\cdot (-ie_0)\subset \dSC$. Note that $\Xi$ and $\oXi$ are disjoint, but
 the de Sitter space lies on the boundary of both domains. 
 
 The same arguments as above can be applied to $\oXi$, leading
 to a family of sesquiholomorphic kernels
 \[\oPsi_\lambda (z,w) ={}_2F_1\Big(\rho+\lambda,\rho-\lambda; \frac{n}{2}; \ltfrac{1+[z,\Bar{w}]}{2}\Big), \quad  z,w \in \oXi .
\] 
 
 These kernels also have distributional boundary values  for fixed $x \in \dS$. We denote the limiting distribution
 by $\oPsi_{\lambda,x}$. It is again analytic on $\Scut$ with a jump singularity on the cut. But it turns out that the wavefront set
 of $\Psi_{\lambda,x}$ and that of $\oPsi_{\lambda,x} $ point in different directions, which in particular implies that $\Psi_{\lambda,x}$ and
 $\oPsi_{\lambda,x}$ are linearly independent.
 
 We now state our main results in more detail. We consider the curve
 \begin{equation}\label{eq:zt}
 z_t = i\cos t e_0 + \sin t e_n \in \Xi, \quad t\in (-\pi/2, \pi/2) 
 \end{equation}
 and note that 
 \[\lim_{t\nearrow \frac{\pi}{2}}z_t =e_n\quad\text{and}\quad \lim_{t\searrow -\frac{\pi}{2}}z_t = -e_n .\]
 Similarly, we set
  \begin{equation}\label{eq:barzt}
 \oz_t = - i\cos t e_0 + \sin t e_n =\overline{z_t}\in \oXi, \quad t\in (-\pi/2, \pi/2) 
 \end{equation}
 and note that  we again have
 \[\lim_{t\nearrow \frac{\pi}{2}}\oz_t =e_n\quad\text{and}\quad \lim_{t\searrow -\frac{\pi}{2}}z_t = -e_n .\]

 We write, as mentioned earlier, $v=(v_0,\bv)$ for $v\in \R{n+1}$ with $v_0\in\R{}$ and $\bv\in\R{n}$.
 We identify $\rT^*(\dS )$ with $\rT(\dS)$ and consider the wavefront set as a subset of 
 \[\rT(\dS) =\{(x,v)\in \R{n+1}\times\R{n+1}\:
 x\in \dS \text{ and } [x,v]=0\}.\]

\begin{manualtheoremx}{1}\label{thm:1} Let $\rho =(n-1)/2$ and  $G = \SO_{1,n}(\R{})_e$. Let $z_t = i\cos t e_0 + \sin t e_n \in \Xi$ and $\overline{z_t}\in \oXi$ for $t\in (-\pi/2, \pi/2)$. Then the following holds
	for $\lambda \in \C{} \setminus\{ (\rho + \mathbb{N}_0) \cup (-\rho - \mathbb{N}_0)\}\colon$
	\begin{enumerate}[leftmargin=*,label={\arabic*)}]
		\item  For $x\in \dS$,  let $g_1,g_2 \in G$ be such that $x = g_1e_n=g_2e_n$. Then
		\[\Psi_{\lambda,x} = \lim_{t \nearrow \pi/2}\Psi_\lambda(g_1z_t,\cdot)=\lim_{t \nearrow \pi/2}\Psi_\lambda(g_2z_t,\cdot),\]
		and 
		\[\oPsi_{\lambda,x} = \lim_{t \nearrow \pi/2}\oPsi_\lambda(g_1\oz_t,\cdot)=\lim_{t \nearrow \pi/2}\oPsi_\lambda(g_2\oz_t,\cdot)
		\]
		exist in the space of distributions. In particular, these limits are well defined for $x \in \dS$. 
		
		\item The analytic wavefront sets of the distributions $\Psi_{\lambda,x}$ and $\oPsi_{\lambda,x}$ are given by 
		\begin{align*}
			WF_A(\Psi_{\lambda,x}) &= \{(x,v)\: v_0 < 0,\, [x,v]= [v,v]=0 \}\\
			& \bigcup \{(x +v, (-v_0, \bv) ) \: v_0<0,\,  [x,v]= [v,v]=0\}\\
			&\bigcup \{(x+v, (v_0, -\bv))\: v_0 <0, \,  [x,v]=[v,v]=0\} 
		\end{align*}
		and
		\begin{align*}
			WF_A(\oPsi_{\lambda,x}) &=
			\{(x,v)\: v_0 > 0,\, [x,v]= [v,v]=0 \}\\
			& \bigcup \{(x +v, (v_0, -\bv) ) \: v_0>0,\,  [x,v]= [v,v]=0\}\\
			&\bigcup \{(x+v, (-v_0, \bv))\: v_0 >0, \,  [x,v]=[v,v]=0\}   . 
		\end{align*} 
	\end{enumerate}
	
\end{manualtheoremx}

Here, $\N_0$ is the set of nonnegative integers. Denote by $\DdS$ the space of distributions on $\dS$. If $L\subset \rO_{1,n}(\R{})$ and $\lambda\in\C{}$, then $\DdS^L$ denotes the space of
$L$-invariant distributions, 
\[\DdSl =\{\Psi \in \DdS \: \Delta \Psi = (\rho^2-\lambda^2)\Psi\}\]
and $ \DdSl^L = \DdSl \cap \DdS^L$. A distribution $\Psi$ is {\it $H$-spherical}, respectively {\it $H^\prime$-spherical}, if there exists $\lambda\in\C{}$ such that
$\Psi \in \DdSl^H$, respectively  $\Psi \in \DdSl^{H^\prime}$, where $H$ and $H^\prime$ are as before.  It is known that $\dim \DdSl^H=4$ (see \cite{OSe80}) and $\dim\DdSl^{H^\prime}= 2$ (see \cite[Prop. B.4.1]{D09}).  In this paper, the basis elements are obtained explicitly as distributional boundary values of the sesquiholomorphic kernels $\Psi_\lambda$ and $\oPsi_\lambda$ on the crown domains, providing a geometric construction.  This construction differs from the approaches via invariant differential equations and Poisson transforms in \cite{OSe80} and generalized Gelfand pairs and spherical distributions in \cite{D09}. Moreover, in this paper, the individual basis elements are distinguished by their analytic wavefront sets. The following results were previously announced, without proofs, by the second author in \cite{Si25}.

To simplify the
notation, let $\Psi_\lambda^\pm= \Psi_{\lambda,\pm e_n}$ and
$\oPsi_\lambda^\pm = \oPsi_{\lambda,\pm e_n}$. Using the wavefront sets, we characterize the spherical distributions in the following two   theorems:    
\begin{manualtheoremx}{2}\label{thm:2}  Let $n \geq 2$ and let $\lambda\in \C{} \setminus \{(\rho + \N_0) \cup (-\rho-\N_0)\}$. Then the following holds:
	\begin{enumerate}[leftmargin=*,label={\arabic*)}]
		\item The distributions $\Psi_{\lambda}^+ $, $\Psi_{\lambda}^- $, $\oPsi_{\lambda}^+  $ and $\oPsi_{\lambda}^- $  are spherical 
		and form a basis for the space $\DdSl^H$.
		
		\item The distributions in part (1) 
		% $\Psi_{\lambda, e_n} $, ${\oPsi}_{\lambda, e_n} $, $\Psi_{\lambda, -e_n}$ and $\oPsi_{\lambda, -e_n}$
		are 
		positive definite if and only if $\lambda \in i\R{} \cup (-\rho,\rho)$. 
		
		\item Let $\Theta \in\DdSl^H$, $\Theta\not= 0$. Then
		\[ WF_A(\Theta)\subset WF_A(\Psi_{\lambda}^+), \sqcup WF_A( \oPsi_{\lambda}^+) \sqcup 
		WF_A(\Psi_{\lambda}^- )\sqcup WF_A (\oPsi_{\lambda}^-). \]
		Furthermore, if  $WF_A(\Theta)$ coincides with the analytic wavefront set of any of  the distributions in part (1),  then $\Theta$ is
		a constant multiple of that distribution.
	\end{enumerate}  
\end{manualtheoremx}

\begin{manualtheoremx}{3}\label{thm:3}
	Let $n \geq 2$ and let $\lambda\in \C{} \setminus \{(\rho + \N_0) \cup (-\rho-\N_0)\}$. Then the following holds:
	\begin{enumerate}[leftmargin=*,label={\arabic*)}]
		\item The distributions $\Theta_{\lambda}^+ := \Psi_{\lambda}^+ + \oPsi_{\lambda}^+ $ and 
		$\Theta_{\lambda}^- := \Psi_{\lambda}^- + \oPsi_{\lambda}^-$  are $H^\prime$-spherical and form
		a basis for $\DdSl^{H^\prime} $.
		
		\item The distributions in part  (1) are  positive definite if and only if $\lambda \in i\R{}\cup (-\rho,\rho)$.
		
		\item Let $\Theta \in \DdSl^{H^\prime}$, $\Theta \not=0$. Then
		\[WF_A(\Theta)\subset WF_A(\Theta_{\lambda}^+)  \sqcup WF_A(\Theta_{\lambda}^-).\]
		Furthermore, if
		$WF_A(\Theta)$ coincides with the analytic wavefront set of any of the distributions in part (1), then  $\Theta$ is 
		a constant multiple of that distribution. 
	\end{enumerate}
\end{manualtheoremx}

Similar ideas for the construction of spherical distributions were first proposed  in a  series of articles \cite{ BV96, BV97,BM96, BM04}
in which the authors considered two ``Lorentzian tuboids'' $\cT^{\pm}$ which are holomorphically equivalent to $\Xi$ and $\oXi$, respectively.  The authors showed
that  the distributional solutions of  \eqref{eq:1} can be extended holomorphically to  $\mathcal{T}^\pm$.    
In particular, it was shown in \cite{BV97} that the 
perikernels are holomorphic in a cut domain of the form $\dS \times \dS \setminus \Sigma$, where $\Sigma$ is the 
set of tuples $(x,y)$ with $[x-y,x-y] \leq 0$. Those articles did not discuss the wavefront set of those distributions. 

This article is organized as follows. The geometry of our setting   is discussed in \cref{sec:ds}. Most of this
material is already covered in \cite{NO20}, except that the conjugate domain $\oXi$ was not discussed there.
The ideas and results presented in this article are closely related to
representation theory, a connection that has been partially discussed in \cite{FNO23,FNO24,NO20}. 
We discuss the kernels $\Psi_\lambda$ and $\oPsi_\lambda$ from the analytic and representation-theoretic viewpoints
in \cref{sec:ker}.
We prove that the boundary values $\Psi_{\lambda,x}$ and $\overline{\Psi}_{\lambda,x}$ are eigendistributions of $\Delta$ on the de 
Sitter space in \cref{sec:bv}. We then review the concept of wavefront sets with some examples in 
\cref{sec: wf}.  In \cref{sec:ws}, we study spherical distributions and prove Theorem 2 and Theorem 3. Part (1) of Theorem~\ref{thm:1} will be proved in two stages. The existence of distributional limits is established in Lemma~\ref{lemma:psi}, while in \cref{sec:proof} we show that the limits depend only on $x \in \dS$. The remaining parts of Theorem~\ref{thm:1} will also be proved in \cref{sec:proof}.
 
 %The motivation for \cite{NO20} was the reflection positivity on the sphere
%  related to the resolvent of the Laplacian. This leads naturally to the kernels $\Psi_\lambda$ and it was shown in \cite{NO18} that the functions 
%$x\mapsto \Psi_\lambda (x,i e_0)$  are the spherical functions
 % on the real hyperbolic space $\mathbb{H}^n$ and that all positive definite spherical 
 %functions on $\mathbb{H}^n$ are obtained in this way. This related the analysis to the spherical principal
%  series representations of $G$, their distribution vectors and to the general 
%theory of the crown as presented in \cite{AG90,vdBD88,BD92,GK02a,GK02b,GKO04,KSt04, S25}.
%Our results are independent of the representation theory, but we refer to \cite{FNO23} where the 
%analytic continuation and the construction of the spherical
%distribution vectors are discussed from the point of view of representation theory.

\subsection*{Acknowledgements:} We would like to thank the referees for their careful comments and valuable suggestions.
The second author acknowledges the support of the Institut Henri PoincarA(C) (UAR 839 CNRS-Sorbonne UniversitA(C)) and LabEx CARMIN (ANR-10-LABX-59-01).
  
\subsection*{Notation.}  We will use the following notation throughout the article:
%he notations we will be using are as follows:\\
%\par
%\textbf{Notations :} 
\begin{itemize}
 \item  $[z,w] = -z_0w_0 +z_1w_1 +...+ z_nw_n$ for $z,w \in \C{n+1}$,
      \item $\R{1,n} = (\R{n+1},[\cdot , \cdot])$,
      \item $G = \SO_{1,n}(\R{})_e\subset {G}' = \mr{O}_{1,n}(\R{})$,
      \item  $H = \mr{SO}_{1,n-1}(\R{})_e=G^{e_n} \subset H^\prime  = \rO_{1,n-1}(\R{})=\rO_{1,n}(\R{})^{e_n}$, 
  \item $\rJ = \begin{pmatrix} -\rI_n & 0\\ 0 & 1\end{pmatrix}$ and $\tau : G^\prime \to G^\prime $, $\tau (g) = \rJ g\rJ$.
     \item  ${K} = \mr{SO}_{n}(\R{})\subset G$,
     \item $G_{\C{}} = \mr{SO}_{n+1}(\C{})$,
     \item  $K_{\C{}} = \mr{SO}_{n}(\C{})$,
    
     \item $\dS =\{x\in \R{n+1} : [x,x]=1\} \simeq G/H \simeq G'/H'$,
 \item $ \dSC = \{z\in \C{n+1}: [z,z]=1\} = G_{\C{}}/K_{\C{}}$,
  \item $\mathbb{H}^n = \{ix\in i\R{n+1} : x_0 > 0,  [x,x] =-1 \} \simeq G/K\subset \dSC$,
  \item $\oHn= \{ix\in i\R{n+1} : x_0 < 0,  [x,x] -1\} \simeq G/K \subset \dS_{\C{}}$,
   \item $\Omega = \{v \in \R{n+1}: [v,v]<0, v_0 >0\}$,
    \item $\rT_\Omega = \R{n+1} + i\Omega$,
 % \item $\mathbb{S}^n = \{(ix_0,\mathbf{x}) : x_0^2 + \mathbf{x}^2 = 1\} \subset i\R{}e_0 + \R{n}$,
%  \item $\Sn_{\pm} = \{(ix_0,\mathbf{x}) \in \Sn: \pm x_0>0\}$,
  \item $\Xi = \dSC \cap \rT_\Omega$ and $\oXi = \dSC\cap \rT_{-\Omega} $, 
   \item $\Gamma^{\pm}(x) = \{ y \in \dS :   [y-x,y-x] < 0, \pm y_0 > x_0 \}$ for $x\in\dS$,
     \item $\Gamma(x) = \Gamma^{+}\cup \Gamma^{-}$,
     \item $\mathbb{L}_{n-1} = \{ v \in \R{n}: -v_0^2+v_1^2+...+ v_{n-1}^2=0\}$,
   
    \item $\sigma(z) = \overline{z} $
 \item $\rho = (n-1)/2$ for $n \geq 2$.
 \item $\N_0 = \{0,1,2,\ldots\}$.
\end{itemize} 
%---------------------------------------------------------------------------------------------------------------------------------------------------------------------------------------------------------------------------------------------------------------------------------------------------------------------------------------------------
\section{ The De Sitter Space, the crown, and the real Hyperbolic spaces} \label{sec:ds}
%%\vspace{3mm}

 \noindent
In this section we recall some basic geometric facts concerning the real hyperbolic space $\Hn$, the de Sitter space $\dS$
and the crown of $\Hn$. The material
is well known and we refer to \cite{NO18,NO22,FNO24} for related discussion.    We will often use the notation from the introduction without further comments. We
start by discussing the de Sitter space and then the hyperboloid $\Hn$ and its crown. 
\subsection{The de Sitter space $\dS$}

%We denote by $\R{1,n}$ the space $\R{n+1}$ viewed as a Lorentzian space with the  Lorentzian form $[\cdot, \cdot]$. We say that a vector
 %$v \in \R{1,n}$ is {\cred {\it time-like} if $[v,v]<0$, {\it space-like} if $[v,b]>0$ and   {\it ight-like}   $[v,v] >0$, and $[v,v]=0$.} 
%Let $G' = \mr{O}_{1,n}(\R{})$ be the isometry group of the bilinear form $[\cdot,\cdot]$ and $G=\SO(1,n)_e$ be the connected component of identity of $G'$. Let $\tau : G' \rightarrow G'$ be the involution given by $\tau(g) = JgJ$, where $J$ is the orthogonal reflection in the hyperplane $x_n = 0$. The restriction  $\tau|_{G}$  is also an involution on $G$ which will be denoted by $\tau$ as well.
As in the introduction, we let $\dS =\{x\in\R{n+1}\: [x, x]=1\}$. To understand
$\dS$ as a symmetric and homogeneous space, we recall the groups $G=\SO_{1,n}(\R{})_e\subset G^\prime =\rO_{1,n}(\R{})$. Both groups under the standard linear action on $\R{n+1}$ act transitively on $\dS$  and the stabilizers of $e_n\in\dS$ are given by
 \begin{align*}
    H &= \{h\in G: h e_n = e_n\}\simeq \SO_{1,n-1}(\R{})_e\\
   & = \left\{ \begin{pmatrix} h & 0 \\ 0 & 1\end{pmatrix}\in\SO_{1,n-1}(\R{})_e \: h\in \SO (1,n-1)_e\right\}
\end{align*} 
and
\begin{align*}
   H' &= \{h\in G'\: h e_n = e_n\} \simeq \rO_{1,n-1}(\R{})\\
    &= \left\{ \begin{pmatrix} h & 0 \\ 0 & 1\end{pmatrix}\in \rO_{1,n}(\R{})\: h\in \rO_{1,n-1}(\R{})\right\} .
    \end{align*}
Both groups have the same Lie algebra 
\[\fh =  \fg^\tau = \left\{\begin{pmatrix} X & 0\\ 0 &0\end{pmatrix} \: X\in \so_{1,n-1}(\R{}) \right\}\simeq \so_{1,n-1}(\R{})\]
where $\tau $ is the involution $\tau (X) = JXJ$, $J (\bv,v_n) = (-\bv,v_n)$.
Moreover, $H = (G^{\tau})_e$ is the connected component of the fixed point group of $\tau$ in G and $(G^{\prime \tau})_e\subsetneq 
H'\subsetneq G'^{\tau}$ has index two in  $G'^{\tau}$.  We have
\[\dS =\{x\in\R{n+1}\: [x,x]=1\}\simeq G/H \simeq G^\prime /H^\prime,\]
realizing $\dS$ as a symmetric space. Finally,
with 
\begin{equation}\label{eq:q}
\fq= \left\{\alpha (v) = \begin{pmatrix} 0 & 0 & v_0 \\
0 &  \mathbf{0}_{n-1} & \bv\\ 
v_0&-\bv^\top & 0 \end{pmatrix}\:v= \begin{pmatrix}v_0\\\bv\end{pmatrix}\in\R{}\times \R{n-1} \right\}\simeq \R{n},\end{equation}
we have $\fq= \fg^{-\tau} $
and $\fg = \fh \oplus \fq$.

The tangent space of $\dS$ at the point $x\in\dS$ is given by
\[\rT_x(\dS) =\{(x,v)\: v\in\R{n+1}\text{ and } [x,v]=0\}.\]
In particular,
\[
\rT_{e_n} (\dS) \simeq % \\
   \fq=\fg^{-\tau}=\{X\in \fg\: \tau (X) =-X\}\simeq \R{1,n-1 }.\]
 In the following we identify $\rT_{e_n}(\dS)$ with $\R{1,n-1}$ using the projection $(e_n,v) \mapsto v$.
 
The tangent bundle is  
\[\rT (\dS) = \{(x,v)\in \R{n+1}\oplus \R{n+1}: x\in \dS\text{ and } [x,v]=0\}.\]
It carries a $G^{ x} = \{g\in G \: g x= x\}$ invariant Lorentzian form $\beta_x$ given by
\[\beta_x(v,w) = [v,w],\]
which turns $\dS$ into a Lorentzian manifold. Moreover, $\rT(\dS)$ is a $G$-homogeneous vector bundle with
the $G $ action $g\cdot (x,v) = (g x, g v)$. We have the isomorphisms
\[\rT(\dS) \simeq_{G} G\times_{H} \fq\simeq_{G^\prime} G^{\prime}\times_{H^\prime} \fq .\]

 We write $\ell_g$ for the diffeomorphism $x\mapsto g x$.
 
 \begin{lemma}
     If $(x,v), (y,w)\in \rT(\dS)$ with $[v,v]=[w,w]$, then
there exists $g\in G$ such that $g{ \cdot } (x,v) = (y,w)$.
 \end{lemma}

\begin{proof}  We can assume that $(x,v) = (e_n,v)\in \rT_{e_n}(\dS)\simeq \R{1,n-1}\times \{0\}$. Let $g\in G$ be such
 that $ge_n =y$. Then   $g^{-1}\cdot (y, w) \in \rT_{e_n}(\dS)$.
 Moreover,  $[g^{-1} \cdot w, g^{-1}\cdot w] = [w,w] =[v,v]$ by  assumption.  Since $H=\SO_{1,n-1}(\R{})_e$ acts transitively on the 
 level sets of $[\cdot , \cdot ]$ in $\R{n}$, there exists $h \in H$ such that $h^{-1}\cdot(g^{-1}\cdot w) = v$. It follows that
 \[gh\cdot (e_n,v) = (g\cdot e_n, g\cdot (h\cdot v))=(y,w) .\qedhere\]
\end{proof}
 
As in \cite[(7.22)]{NO18} and \cite[(2.3), p.15]{NO20}, we define the analytic functions $\rC,\rS:\bbC \rightarrow \bbC$ by
$$\rC(z):= \sum_{k=0}^\infty \frac{(-1)^k}{(2k)!}z^k = \cos \sqrt{z} \quad \text{and} \quad  \rS(z):=\sum_{k=0}^\infty \frac{(-1)^k}{(2k+1)!}z^k
=\frac{\sin \sqrt{z}}{\sqrt{z}} .$$
 Note that  both functions are well-defined since the functions $y \mapsto \cos (y)$ and $ \sin (y)/y$ are even.

 \begin{remark}\label{rem:2.2} In the following we assume that $(x,v)\in \rT_x(\dS)$.  We note the following:
 \smallskip
 
 \noindent
 1) Using that $[x,v]=0$ and writing $z=\sqrt{[v,v]}$,
 \begin{align*}
 [\rC(z)x + \rS(z)v, \rC(z)x + \rS(z)v]
 &=[\rC(z)x,\rC(z)x ]+ [\rS(z)v,  \rS(z)v]\\
 & = \cos^2(z)[x,x] +\frac{\sin^2(z)}{[v,v]}\, [v,v]\\
 &=\cos^2(z) +\sin^2(z)=1 .
 \end{align*}
Hence $\rC(z)x + \rS(z)v\in\dS$.
\smallskip

\noindent
2)  If $[v,v] =0$, then $\rC([v,v])x+\rS([v,v])v = x+v$. For
the case $[v,v]\not=0$, for simplicity, take $x=e_n$ and $v =(r\bv ,0)$ with $[\bv,\bv]=\pm 1$ and $r>0$.  If $[\bv,\bv]=1$, then
\[\rC([v,v])x+\rS([v,v])v = (\sin (r) \bv, \cos (r))\]
is a circle. On the other hand, if  $[\bv,\bv]= -1$, then 
\[\rC([v,v])x+\rS([v,v])v = (\sinh (r) \bv ,\cosh (r))\]
is a hyperbola and the map is injective as a function of $r\in (0,\infty)$.
\end{remark}

 In the following we simplify our notation by writing $\Exp_x(x,v) = \Exp_x(v)$ for $(x,v) \in \rT_x(\dS)$.
\begin{lemma}[N{\'O}20]\label{lemm:exp1} The exponential map  $\Exp_x : \rT_x(\dS)\to \dS$ is given by 
$$\mathrm{Exp}_x(v) = C([v,v])x + S([v,v])v, \quad v \in \rT_x (\dS), $$
and we have
\begin{equation}\label{eq:ellg}
\ell_g \circ \Exp_x = \Exp_{g x} \circ (d\ell_g)_x.
\end{equation}
\end{lemma}  

\begin{proof} This is well known, but here is the main idea of the proof. Using $[g\cdot v,g\cdot v]= [v,v]$ and the definition of the action of 
	$G$ on $\rT(\dS)$, we obtain \eqref{eq:ellg}. We can therefore assume that $x=e_n$. Then a simple calculation shows that
\[\exp \alpha(v) = \left( \begin{array}{c|c}
    \begin{array}{cc}
        \rC([v,v] ) & 0  \\
        0 & \rI_{n-1}
    \end{array} & \rS([v,v])v \\
    \hline \\[-2ex]
    -\rS([v,v]) v^T & \rC([v,v])
\end{array} \right)\]
and hence
\[\Exp_{e_n}(\alpha (v)) = \exp(\alpha (v))e_n = \rS ([v,v])v + \rC([v,v])e_n.\qedhere\]
\end{proof}

\begin{lemma}\label{lemma: exp_x}
Let $U_{e_n} =\{v\in \rT_{e_n}(\dS) : |[v,v]|<\pi^2\}$. For $x=g \cdot e_n $,
set $U_x = d\ell_g U_{e_n}$ and $V_x = \mathrm{Exp}_x U_x \subset \dS$. Then the following holds: $V_x$ is open and $\Exp_x : U_x \to V_x$ is a real analytic diffeomorphism.
\end{lemma}
\begin{proof}
Clearly, the map is real analytic. By \eqref{eq:ellg}, we can assume that $x=e_n$. To show that $\Exp_{e_n}$ is
injective on $U_{e_n}$, we use part 2 of Remark \ref{rem:2.2} and Lemma \ref{lemm:exp1}. If $[v,v]=0$, then
$\Exp_{e_n} (\alpha (v))= e_n + v$, which is injective as a function of $v$. Otherwise, we write $v=r\bv$
with $[\bv,\bv]=\pm 1$ and $r>0$. The calculation in Remark \ref{rem:2.2} shows that the 
\[(0,\pi)\times \{\bv \: [\bv,\bv]=\pm 1\}\ni (r,\bv) \mapsto \begin{cases} (\sin (r) \bv ,\cos (r) )& [\bv,\bv]=1\\
(\sinh (r)\bv, \cosh (r))& [\bv,\bv]=-1\end{cases}\]
is injective. 

Identifying $\rT_{e_n}(\dS)$ with $\R{n}$ via $\bv\mapsto (\bv,0)$ and identifying $\rT_{\bv}(\R{n})$ with
$\R{n}$, one easily sees that with $a= \sqrt{[v,v]}$ we have
\[(D \Exp_{e_n})_{\bv}(w) = \frac{-\sin a}{a}[\bv,w]e_n + \frac{[\bv,w](a\cos a - \sin a)}{a^3} \bv +\frac{\sin a}{a}w \]
which is injective as a function of $w$ as long as $v=(\bv,0)\in U_{e_n}$.

\end{proof}

We obtain from Lemmas \ref{lemm:exp1} and \ref{lemma: exp_x} that the family $\{U_x, \mr{Exp}_x\}_{x \in \dS}$ forms a real analytic  atlas
and $G$ clearly acts analytically in those coordinates.

A similar description for the tangent bundle of $\dSC$ holds and we leave it to the reader to formulate
the exact statement.

\subsection{Invariant differential operators on $\dS$} \label{sec:diffoper}
Assume that the Lie group $L$ acts on a  manifold $X$ by analytic diffeomorphisms $ \ell_g(x)=g\cdot x$, $g\in L, x\in X$. We say
that the differential operator $D$ on $X$ is $L$-{\it invariant} if for all $f \in C_c^\infty (X)$ and
all $g\in L$,
we have $D(f \circ \ell_g) = (Df)\circ \ell_g$. If $D_1,D_2$ are $L$-invariant, then so is $D_1D_2$. We denote by 
$\D (X)$ the algebra of invariant differential operators. It is known that $\D (\dS) = \C{} [\Delta_1]$, the algebra of polynomials in the
Laplacian-Beltrami operator $\Delta_1$ (see \cite{F79,D09}). The generator of the algebra $\D (\dS)$ can also be
defined in the following way.

First let 
\[\square_{n+1} = -\dfrac{\partial^2}{\partial x_0^2} + \sum_{j=1}^n \dfrac{\partial^2}{\partial x_j^2}\]
in $\R{1,n}$.
Let $\varphi \in C_c^\infty (\R{})$, $\varphi = 1 $ in a neighborhood of $1$, and $\varphi (t) =0$ for
$|t-1|> 1/2$. For $f \in C_c^\infty (\dS)$, define
\[\tilde f (x) = \varphi ([x,x]) f(x/| [x,x]|^{1/2}),\quad x\in \R{1,n}.\]
Then $\tilde f \in C^\infty_c(\R{1,n})$. Define
\[\Delta_2 f :=  (\square_{n+1}\tilde f)|_{\dS}. \]
This is a well-defined  second-order  $G',G$-invariant differential operator on $\dS$ (see \cite[p. 110, 160]{D09}). 

The third realization of the generator is to use the coordinates from Lemma \ref{lemma: exp_x}. 
Note that $\square_n$ is a well-defined $H$, $H'$-invariant differential operator
on $\rT_{e_n}(\dS) \simeq \R{1,n-1}$. Define  
\[(\Delta_3 f)\circ \Exp_{e_n} := \square_{n}(f\circ \Exp_{e_n}),\quad f\in C_c^\infty (V_{e_n}).\]
As $\square_{n}$ is $H$, $H'$-invariant, it is well defined on $U_{e_n}$. For $f \in C_c^\infty(V_x)$ and $g\in G^\prime$ such
that $x=ge_n$ , define 
\[(\Delta_3 f)\circ \mr{Exp}_{x} = (\square_n(f \circ \ell_g \circ \mr{Exp}_{e_n}))\circ d\ell_g^{-1}. \]
By applying a partition of unity subordinated to the open
cover $\{V_x\}_{x\in \dS}$, we have a well-defined second-order $G'$, $G$-invariant differential operator  $\Delta_3$ on $\dS$.

 Since all three operators are second-order, annihilate constants and preserve positivity, they agree up to a multiple of a strictly positive constant and hence have the same sets of eigenfunctions.
 From now on we will simply denote a choice of generator by   $\Delta$. 
 
\subsection{The Crowns $\Xi$ and ${\oXi}$}

The complex crown  $\Xi$ of a Riemannian symmetric space $X=G/K$, with $G$ linear and $K$ maximal compact, is a natural $G$-invariant complex open domain 
in the complexification $G_\bbC/K_\bbC$ with the property that eigenfunctions of the algebra of
invariant differential operators extend to $\Xi$. It was first introduced in \cite{AG90} and then  studied by several
authors, including \cite{GK02a,GK02b,KSt04}, whose work is particularly relevant for us. In particular, it was shown in \cite{GK02b} 
that every non-compactly causal symmetric space (see \cite{O91,HO97} for the definition and more information), including the 
de Sitter space, can be
realized as an open orbit in the boundary of the crown. The crown appeared again naturally in \cite{NO20} in 
relation to reflection positivity. We will collect results from \cite{NO20}  here. Our discussion, as well as the material in
\cite{NO20}, is independent of the general theory and builds on the concrete setting considered here.
  
The Riemannian symmetric space associated to $G=\SO_{1,n}(\R{})_e\subset G^\prime = \rO_{1,n}(\R{})$ is
\[\mathbb{H}^n = G\cdot ie_0\simeq G/K\]
with
\begin{align*}
K &=\SO_n (\R{})_e\simeq \left\{\begin{pmatrix} 1& 0 \\ 0 & B\end{pmatrix} \: B \in \SO_{n}(\R{})_e\right\} %\subset\\
%K^\prime &= \rO_n (\R{})_e\simeq \left\{\begin{pmatrix} 1& 0 \\ 0 & A\end{pmatrix} \: A \in \rO_n(\R{})_o\right\} 
\end{align*}
which is the stabilizer of $ie_0$ in $G$. Note that $\SO_n(\R{})$ is connected for $n\ge 2$. The group $K$ is the fixed point group of the
Cartan involution 
\[\theta (g) :=  (g^{-1})^\top\]
defined on $\mr{SO}_{1,n}(\R{})_e$.

Let $h =E_{0n} + E_{n0}\in \fg= \so (1,n)$ be the linear map
\[h(x_0,x_1, \ldots , x_{n-1}, x_n) = (x_n,0, \ldots ,0,x_0).\]
Then $\ad h, X\mapsto [h,X]$, has the eigenvalues $0, 1,-1$, leading to the eigenspace decomposition of $\fg$:
\begin{equation}\label{eq:g1}
    \fg=\fg_{-1}\oplus \fg_0 \oplus \fg_{+1},
\end{equation}
 where
\[\fg_j = \{X\in \fg \: \ad h(X) =  jX\},\qquad j=0,1,-1.\]
As $[\fg_i,\fg_j]\subset \fg_{i+j}$, it follows that each of the spaces $\fg_{\pm 1}$ is $\fg_0$-invariant.

The crown of the Riemannian  hyperbolic space $\Hn =G\cdot (ie_0)$  is defined to be
\[\Xi = G\exp (i(-{\pi}/{2},{\pi}/{2}) h )  \cdot ie_0 = G\cdot \{(i\cos t, 0,\ldots , \sin t):
|t|<\pi /2\}. \]
Similarly for the conjugate space $\oHn$, the crown is given by
\[\oXi = G (\exp i(-\pi/2, \pi /2)h  )\cdot( -ie_0)  = \sigma (\Xi),\]
where $\sigma$ is the complex conjugation on $\C{n}$, i.e., $\sigma(z) = \overline{z}$.

We write $z_t = \exp (-ith) \cdot ie_0=i\cos (t)e_0 + \sin (t)e_n$ and similarly
$\oz_t= \exp (ith) \cdot (-ie_0) = -i\cos (t)e_0 + \sin (t)e_n$, $t\in  (-\pi/2, \pi/2)$. 
Then for   $z \in \Xi$ and $w \in \oXi$, there exists $g,g' \in G$ such that 
\begin{equation}\label{eq:zw}
    z = g z_t \quad \text{and} \quad w = g' \overline{z_t}.
\end{equation}
 
 We now recall the description of $\Xi $ and its properties; see \cite{NO20}. The corresponding statements
for $\oXi$ follow by taking the complex conjugation  $\sigma$.

Consider the open future light cone $\Omega$ given by
 $$\Omega = \{x \in \R{1,n}: [x,x]<0, x_0>0\}. $$
The corresponding future and past tube domains in $\mathbb{C}^{n+1}$ are given by
 
 $$\rT_{\Omega} = \R{1,n} + i\Omega  \quad \text{and}\quad \rT_{-\Omega} = \R{1,n} - i\Omega .$$

 Realize the $n$-sphere by
 $\Sn =\{x\in i\R{}e_0 +  \{0\}\times \R{n}: [x,x]=1\}$ and write
 \[\Sn_+ = \{ x\in \Sn: x_0>0\}\quad{and} \quad  \Sn_- = \{ x\in \Sn: x_0 < 0\} = \sigma (\Sn_+).\]

\begin{lemma}[N{\'O}20] \label{lemma:crown}
The crowns can be described as
 \begin{align*}
 \Xi & =  \{ u + iv : [u,u]-[v,v]=1, [u,v]=0,[v,v]<0, v_0>0\};\\
 &= \rT_\Omega \cap \dS_{\C{}} = \rT_\Omega \cap \Sn_{+,\C{}}=G \Sn_+ 
 \\
 \oXi  &= \{ u - iv :[u,u]-[v,v]=1, [u,v]=0,[v,v]<0, v_0 > 0\}\\
 &= \rT_{-\Omega} \cap \Sn_{-,\C{}} \cap \dS_{\C{}}=G \Sn_-=\rT_{-\Omega}. 
\end{align*}
\end{lemma}
\begin{proof} The first part is \cite[Lem. 3.1]{NO20} and \cite[Prop. 3.2]{NO20}. The second part follows by applying
$\sigma$ to $\Xi$.
\end{proof}

%The following is immediate from above lemmas. 
%\begin{corollary}  
%The crown domains $\Xi$ and $\oXi$ are open subsets of $\Sn_{\C{}}$, hence they are complex
% manifolds. The group $G$ acts on $\Xi$ as holomorphic maps and acts on $\oXi$ as anti-%holomorphic maps.
%\end{corollary}

The following proposition and lemma  are the key tools for extending the
kernels   $\Psi_{\lambda}$ and $\oPsi_{\lambda}$ to the crowns $\Xi$ and  $\oXi$, respectively, and
then taking the boundary values (see \cref{sec:ker}).
\begin{proposition}[N{\'O}20] We have
 \[ \{[z,\sigma(w)]: z,w \in \Xi \} = \C{} \setminus [1,\infty)  =\{[z,\sigma (w ) ] : z,w \in  \oXi \} .\]
\end{proposition}
\begin{proof} The crown is invariant under the conjugation $z\mapsto -\sigma (z)$. The claim therefore
follows from \cite[Lem. 3.5]{NO20}, since the Lorentz form in \cite{NO20} is the negative
of the form considered here. The claim for $\oXi$ follows from the first part since
$\oXi = \sigma (\Xi)$.
%
 %   The proof of first equality is in \cite{NO18}. Now, 
 %$[\tz,\bar{\ti{w}}] = \overline{[\Bar{\tz}, \ti{w}]}$. Since 
 %$\Bar{\tz} \in \Xi$, from the second equality follows from the first equality.
\end{proof}

For $U\subset \bbC^{n+1}$, denote by $\cl{U}$ the closure of $U$ in $\bbC^n$. The boundary
is $\partial U = \cl{U} \setminus U$.

\begin{lemma}\label{lemma:crownboundary} 
    The boundaries of $\Xi$ and $\oXi$ in $\dSC$ are given by
  \begin{align*}
      \partial\,\Xi &= \{x+iy: x,y \in \R{1,n}, \,[x,x]=1, [y,y]=0,y_0 \geq 0, [x,y] =0\};  \\
       \partial\, \oXi &= \{x-iy: x,y \in \R{1,n}, \,[x,x]=1, [y,y]=0,y_0 \geq 0, [x,y] =0\}\\
       &=\sigma (\partial \Xi).
  \end{align*}  
\end{lemma}

\begin{proof} The first part is \cite[Lem. 3.7]{NO20} and the second claim then follows from
$\oXi = \sigma (\Xi)$.
%Suppose that $y_0 = 0$, then $-y_0^2 + \mathbf{y}^2 \leq 0$ implies that $\mathbf{y}^2 \leq 0 $, hence $y = 0$. Let $x \in \mr{RHS}$.  
%Then there exist $\Tilde{y}$ such that $[x,\Tilde{y}] =0$, $[\Tilde{y},\Tilde{y}]=-1.$  Let $z_\epsilon = ( \sqrt{1+\epsilon^2})^{-1}({x + i\epsilon\Tilde{y}})$ lies in $\Xi$. As $\epsilon \rightarrow 0$, $z_\epsilon \rightarrow x$. 
%\par
%Let $z=x+iy$. Suppose $y_0>0$, $[y,y]=0$ and $[x,x]=1$. We can make $v_0 = 0$ by acting some element of $G$. Define $z_\epsilon = (\sqrt{1+2\epsilon y_0 + \epsilon^2})^{-1}(z+ \epsilon e_0) \in \Xi$ and $z_\epsilon \rightarrow Z$ as $\epsilon \rightarrow 0$.     The proof of the second one is also the same.
\end{proof}

\begin{corollary} $\dS=\partial \, \Xi \cap \partial\, \oXi$.
\end{corollary}
\begin{proof}  The above description of the boundary implies,  by taking $y=0$, that
$\dS \subset \partial\, \Xi \cap \partial\,  \oXi$. For the other inclusion, we get from Lemma~\ref{lemma:crownboundary} that $z= x+iy\in \partial\,  \Xi \cap \partial\,  \oXi$ if and only if $y_0 = 0$. Hence, 
$0 = [y,y]$   leads to $y_0^2 =\|\by\|^2=0$, so $y=0$. The claim now follows from   $[x,x]=1$. 
\end{proof}

%------------------------------------------------------------------------------------------------------------------------------------------------------------------

%----------------------------------------------------------------------------------------------------------------------------------------------------------------
%%\vspace{5mm}
\section{A representation-theoretic overview }\label{sec:ker}

 We have already mentioned that some of the ideas in this article are closely related to
representation theory. In \cite{NO20}, the kernels $\Psi_\lambda$ were introduced, arising 
from the application of representation-theoretic ideas to constructive quantum field theory in the form of
reflection positivity. These kernels are positive-definite for parameter $\lambda \in i\R{}\cup (-\rho,\rho)$  and  $\varphi_\lambda  = \Psi_\lambda (\cdot ,ie_0)$ is the spherical function on $\Hn$ with spectral
parameter $\lambda$. 
Thus, in this section we will only consider $\lambda \in i\R{}\cup (-\rho,\rho)$.

The existence of the  limit of spherical
functions in the space of distributions  on non-compactly causal symmetric spaces was already proved  in \cite{GKO03,GKO04} using 
advanced tools from representation
theory, in particular the {\it automatic continuation theorem}  of van den Ban, Brylinski and Delorme; see
\cite[Thm. 2.1]{vdBD88} and \cite[Thm. 1]{BD92}, also restated in \cite{NO18} for the specific case of $\dS$.
A different and less abstract proof was given in \cite{FNO23} and a more general version in \cite{S25}.  
But none of those articles considered the detailed analysis of the singularities carried out in this article. Our method is
also different from those articles.

\subsection{The positive definite kernels $\Psi_\lambda$ and $\oPsi_\lambda$}

We recall that reflection positivity on the sphere \cite{NO20} (see also \cite{NO22}) led to the
positive definite kernel $\Psi_\lambda$ for $\lambda \in i[0,\infty) \cup [0,\frac{n-1}{2})$, introduced in \eqref{eq:1.1}
(also denoted by $\Psi_m$ and normalized differently in \cite[Thm. 4.12, Thm. 4.19]{NO20}). For $z,w \in \Xi$, this kernel is given by
\begin{equation}\label{def:PsiLambda}
\Psi_\lambda (z,w) = {}_2F_1\left( \rho + \lambda, \rho-\lambda ; \frac{n}{2}; 
\frac{1+[z,\sigma (w)]}{2}\right), \quad \rho = \frac{n-1}{2}.
\end{equation}
As both $n$ and $\lambda$ are fixed most of the time, we simplify our notation and write
\[\Fl (z) = \Fl (n,\lambda ;z) = {}_2F_1\left( \frac{n-1}{2}+ \lambda,  \frac{n-1}{2}-\lambda ; \frac{n}{2}; z\right).\]
Here ${}_2F_1(a,b;c;z)$ denotes the Gauss hypergeometric function
\[{}_2F_1(a,b;c;z) = \sum_{n=0}^\infty \frac{(a)_n(b)_n}{(c)_n}\frac{z^n}{n!},\]
where $(d)_n = d(d+1)\cdot(d+n-1)$,  $c\not\in -\N_0$ and $|z|<1$. Recall the following fact:  
\begin{theorem} The hypergeometric function ${}_2F_1$ extends to a holomorphic function
on $\bbC\setminus [1,\infty)$. Furthermore, for a fixed $z\in \bbC\setminus [1,\infty)$, the
function
\[(a,b,c) \mapsto \frac{{}_2F_1(a,b;c;z)}{\Gamma (c)}\]
is an entire function of $a$, $b$ and $c$, where $\Gamma(c)$ is the Gamma function.
\end{theorem}
\begin{proof} See \cite[p.241,245--246]{LS}.
\end{proof}

We note that by G-invariance $\Psi_\lambda $ is uniquely determined by 
\[\varphi_\lambda (x) =
\Psi_\lambda (x,ie_0)=  {}_2F_1 \left( \rho + \lambda, \rho -\lambda ; \frac{n}{2}; 
\frac{1+ix_0}{2}\right),\quad ix\in\Hn,\]
which is the positive definite spherical function on $\Hn$ with spectral parameter $\lambda$.  

%Using the fact that $-\sinh^2 (t/2) = \frac{1-\cosh (t)}{2}$, we have
%\begin{equation}\label{eq:sphf}
%\phi_\lambda( \exp (th) ie_0) = {}_2F_1 \left( \rho + \lambda,  \rho-\lambda ; \frac{n}{2}; 
%-\sinh^2(t/2) \right). 
%\end{equation}

Similarly, we have that the kernel $\oPsi_\lambda$ on $\oXi\times \oXi$ is given by  
$$\oPsi_{\lambda}(z,w) = \, {}_2F_1\Big(\rho+\lambda,\rho-\lambda; \frac{n}{2}; \ltfrac{1+[{z}, \sigma(w)]}{2}\Big), \quad  z,w \in \oXi.$$
This is
uniquely determined by $\oPsi_\lambda(\cdot , -ie_0)|_{\oHn}$, the spherical function on $\oHn$ with spectral parameter $\lambda \in i\R{} \cup (-\rho,\rho)$.   

We collect the main facts in the following theorem. Part (1) is in \cite[Sec. 4]{NO20} and \cite[p.11]{NO22}.

\begin{theorem}\label{thm:kernel} Let $\rho = \frac{n-1}{2}$ and $\lambda\in i\R{}\cup (-\rho,\rho)$. Then
\begin{itemize}
    \item[\rm (1)] The kernel $\Psi_{\lambda}(z,w)$  is a positive definite 
    $G$-invariant kernel  on $\Xi \times \Xi$, holomorphic in the first 
variable and anti-holomorphic in the second variable. 
\item[\rm (2)] Let $z,w\in\Xi$. Then
$\overline{\Psi_\lambda (z,w)} = \Psi_\lambda( \sigma(z),\sigma (w))$. 

    \item[\rm (3)]  The kernel $\oPsi_{\lambda}$ is
a $G$-invariant positive definite kernel on
 $\oXi \times \oXi$, holomorphic in the first variable and
 anti-holomorphic in the second variable. 
\end{itemize}
  
\end{theorem}
\begin{proof} We use the simplified notation ${}_2F_1(u)= {}_2F_1(\rho + \lambda, \rho -\lambda ;n/2; u)$.
(1) is \cite[Thm. 4.12]{NO20} and (2) follows from (1) and the fact that for $\lambda \in i\bbR \cup \bbR$ and $z\in \C{}\setminus [1,\infty)$,
we have
\[\overline{{}_2F_1(z)} =
{}_2F_1(\bar z),\]
as ${}_2F_1(a,b;c;z) = {}_2F_1 (b,a; c; z)$. 
The last claim follows from $\oXi = \sigma (\Xi)$ and the relation
  ${\Psi}_{\lambda}(z,w) = \oPsi (\ow ,\oz)$. Thus, (3) follows from (1) and the arguments given above. 
\end{proof}

\subsection{The boundary value distribution}\label{seBoundVal}  If $L$ is a locally compact topological group and $M $ is a
closed subgroup, then there exists an $L$-invariant measure $\mu_X$ on $X=L/M$ if both $L$ and $M$ are
unimodular. If $L$ is compact, then we normalize $\mu_X$ so that $\mu_X(X) =1$. In particular, there exists
a $G$-invariant measure $\mu_{\dS}$ on $\dS$ (see \cite[Thm. 8.1.1]{S84}) and a normalized $\mr{SO}_{n}(\R{})$-invariant
measure $\mu_{\rS^{n-1}}$. Note that $G$ acts on $\rS^{n-1}$, but the measure $\mu_{\rS^{n-1}}$ is not $G$-invariant. Also, the $G'$-invariant measure on $\dS$ is given in \cite[p. 159]{D09}.  We simplify the notation and write $dy = d\mu_{\dS}(y)$.

Let $z_t = \exp{ ( -  ith) ie_0} = i\cos t e_0 + \sin t e_n$. If $|t| < \pi/2$, then the kernel $\Psi_\lambda(\exp{ (ith) ie_0}, w)$ is 
anti-holomorphic in $w$ and continuous on $\cl{\oXi} \supset \dS$.  Therefore, it defines a distribution
on  $\dS$ given by
\begin{equation}\label{eq:analdist}
    \Theta_\lambda(z_t)(\varphi) = \int_{\dS}\overline{\varphi(y)} \Psi_\lambda( \exp{ (ith) ie_0},y)dy, \qquad \varphi \in C_c^\infty(\dS).
\end{equation}  

 We need some more information to discuss
the connection to representation theory.
Let $\fa =\R \, h$, $\fm = \fz_\fk (\fa)$ and $\fn = \fg_1$, where $\fg_1$ is the $+1$ eigenspace of $\ad h$ (see (\ref{eq:g1})). Then
$\fp = \fm \oplus \fa \oplus \fn$ is a minimal parabolic subalgebra. The corresponding
minimal parabolic subgroup is $P = N_G (\fp) = MAN$, where 
$A = \exp \R{} h$, $N =\exp \fn$, and 
\[M = \left\{\begin{pmatrix}
    1 & 0 & 0\\
    0 & B & 0\\
    0 & 0 & 1
\end{pmatrix}: B \in SO(n-1) \right\}. \] Note that $N$ is abelian in this case. Furthermore,
$G = KAN \simeq K \times A\times N$. Accordingly, we write $g = k(g) a(g) n(g)$. We have
$G/P = K/M = \rS^{n-1}$ and  $G$ acts on $\rS^{n-1}$ by 
$g v= k(g)v$.  

Denote by $\pi_\lambda$ the principal series representation with spectral
parameter $\lambda \in \C{} $. It acts   on the Hilbert space $\mathcal{H}_\lambda=L^2(\rS^{n-1})$ (with the $K$-invariant
probability measure) by
\[\pi_{\lambda }(x) f(v) = a(x^{-1}k)^{-\lambda - \rho} f(x^{-1} v) . \] 
 The representation $\pi_\lambda$ is unitary for $\lambda \in i\R{}$. Furthermore, the representations $\pi_\lambda$ and $\pi_{-\lambda}$ are equivalent.

 Let $v= k(ie_0)\in\rS^{n-1}$,  $x\in G$  and  $f\in L^2(\rS^{n-1})$. Let $\langle\cdot,\cdot \rangle_{L^2}$ be the $L^2$ inner product.
The constant function $e_\lambda (v) = 1$ is $K$-invariant with norm $1$ and
the associated spherical function is
\[\varphi_\lambda (x) = \ip{\pi_\lambda (x)e_\lambda}{e_\lambda}_{L^2} =
\int_{\rS^{n-1}} a(x^{-1}v)^{-\lambda - \rho}dv .\] 
We note that $x\mapsto \pi_\lambda (x)e_\lambda$ is right $K$-invariant. Hence,
$\pi_\lambda (z)e_\lambda$ is well defined  on $\Hn$ and extends to a vector-valued
holomorphic function on $\Xi$. Define a G-invariant kernel $\Psi^1_\lambda$ on $\Xi\times \Xi$ by
\begin{equation}\label{eq:Psi1}
\Psi^1_\lambda (z,w) = \ip{\pi_\lambda (z)e_\lambda}{\pi_\lambda (w)e_\lambda}_{L^2}.
\end{equation}

It was proved in \cite[Sec. 5.1]{NO20} that the kernels $\Psi^1_\lambda$ and $\Psi_\lambda$ coincide on $\Xi\times \Xi$. 
\begin{Remark}
   For each $\lambda \in (-\rho, \rho)$, there exists a different $G$-invariant inner product $\ip{\cdot}{\cdot}_\lambda$
  on $\mathcal{H}_\lambda$. These give rise to the  {\it complementary series representations}.  The kernel function $\Psi^1_\lambda$
  is still defined by \eqref{eq:Psi1}, with $\ip{\cdot } {\cdot }_{L^2}$ replaced by $\ip{\cdot } {\cdot }_{\lambda}$. 
  Again, $\Psi^1_\lambda$ is a $G$-invariant kernel that coincides with $\Psi_\lambda$ (see \cite[Sec. 5.1]{NO20}) on $\Xi \times \Xi$.
\end{Remark}  
For $g \in G$, we have $g\mr{exp}(-ith){ i}e_0 \in \Xi$. The orbit map
\[t \mapsto \pi_\lambda (g\exp (-ith )){ i}e_\lambda, \qquad t \in (-\pi/2,\pi/2), 
\]
is analytic and 
\begin{equation}\label{eq:eH}
e_\lambda^H = \lim_{t \nearrow \pi/2} \pi_\lambda (\exp (-ith )){ i}e_\lambda
\end{equation}
exists in $\mathcal{H}_\lambda^{-\infty}$ and defines an $H$-invariant distribution
vector \cite[Sec. 5]{FNO23}. Furthermore, $\pi^{-\infty}_\lambda (\varphi)e_\lambda^H\in \mathcal{H}^\infty_\lambda$ for
$\varphi \in C_c^\infty (G/H)$; see \cite[Chap. 7]{NO18}.
 Hence,
\[\Theta_\lambda (\varphi ) = \ip{e_\lambda^H}{\pi_\lambda^{-\infty}(\varphi)e_\lambda^H}\]
defines an $H$-invariant distribution. Moreover,
\[\Delta \Theta_\lambda = ( \rho^2 - \lambda^2)\Theta_\lambda .\]

This can be reformulated in terms of the kernel $\Psi_\lambda$. For this, let $z\in \Xi$ and $g\in G$.  Then
\[t\mapsto \Psi_\lambda (z,g\exp (- ith)ie_0)=\Psi_\lambda (z,g(i\cos t e_0 + \sin t e_n))\]
is analytic on an open interval containing $(-\pi/2,\pi/2)$ and has the limit
\[{}_2F_1\left(\rho + \lambda , \rho-\lambda ; \frac{n}{2}; \frac{1 +[z,ge_n]}{2}\right) = \Psi_\lambda (z,ge_n).\]
We then get a
distribution on $\dS$ by
\[\Theta_\lambda(z; \varphi ) = \int_{\dS} \overline{\varphi (y)} \Psi_{\lambda}(z,y)d\mu_{\dS}(y)
= \ip{\pi_\lambda(z)e_\lambda}{\pi^{-\infty}_\lambda (\varphi)e_\lambda^H},\]
where $\mu_{\dS}$ is a $G$-invariant measure on $\dS$.
Taking the limit $z\to e_n$ leads to the eigendistribution $\Theta_\lambda$:
\begin{equation}\label{eq:theta limit}
    \Theta_\lambda (\varphi) = \lim_{t \nearrow \pi/2}  \int_{G/H} \overline{\varphi (y)} 
\Psi_{\lambda}(\exp (-ith)e_n,y)d\mu_{\dS}(y)
\end{equation}

or
\[\Theta_\lambda = \lim_{t \nearrow \pi/2} \Psi_\lambda (\exp (-ith ) e_n, \cdot ) .\]

%---------------------------------------------------------------------------------------------------------------------------------------------------------------------------------------------------------------------------------------------------------------------------------------------------------------------------------------------------}
%\vspace{5mm}
\section{Distributions as Boundary values of holomorphic functions}\label{sec:bv}

In this section we study the kernels $\Psi_\lambda(z,.)$ and ${\oPsi}_\lambda(\bar{z},.)$ for $\lambda \in \C{}$. We show in particular  
that the boundary values of the kernels  $\Psi_\lambda(g z_t ,.)$ and $\oPsi_\lambda(g  \bar{z_t},.)$ are distributions on $\dS$ as $t \to \pi/2$ from below without using the abstract Automatic Continuation Theorem or the existence of the $H$-invariant distribution vector $e^H_\lambda$. 

We use the  notation $\mathcal{D}(\dS) = C_c^{\infty}(\dS)$, $\mathcal{E}(\dS)= C^\infty(\dS)$, both endowed with the
standard topology. The conjugate linear duals  $\cD^* (\dS)$ and $\mathcal{E}^*(\dS)$ are the spaces of distributions and compactly supported distributions on $\dS$, respectively.  

 %-------------------------------------------------------------------------------------------------------------
 \subsection{Kernels $\Psi_\lambda$ and $\oPsi_\lambda$}\label{sec:behaviour}

The following is \cite[Lem. 6.4]{NO20}:
\begin{lemma}\label{l1} We have
\[[\dS ,\Xi]\cap \bbR = [\dS , \oXi]
\cap \bbR = (-1,1).\] 
\end{lemma}

This implies that for $z\in\Xi, w \in \oXi$ and $y\in \dS$, we have
\[\frac{1+[z,y]}{2} \quad \text{and}\quad \frac{1+[w,y]}{2} \not\in [1,\infty).\]
Hence the functions
\[(z,y)\mapsto \Psi_\lambda (z,y), \quad (w,y)\mapsto \oPsi_\lambda (w,y)\]
are analytic.
From this and the properties of the hypergeometric function we get:

\begin{lemma}%\label{cont}
 The kernel $\Psi_{\lambda}$ can be extended continuously to $\Xi \times (\dS\cup\Xi) $ and the kernel $
 \oPsi_{\lambda}$ can be extended continuously to $\oXi \times( \dS\cup\oXi)$ for any $\lambda\in \C{}$.
\end{lemma}
%\begin{proof}
% The proof follows from the properties of hypergeometric function  and \cref{l1}.
%\end{proof}
%\par
%%\vspace{5mm}

%We now note $\bbC^n\setminus [1,\infty)$ is open, so that if $x,y\in \dS$ and $(1+[x,y])/2\not\in [1,\infty)$
%so that $\Psi_\lambda(x,y)$ and $\oPsi_\lambda (x,y)$ are defined, then there
%exists an open set $U$ containing $x$ such that
%\[U\to \bbC, z\mapsto \Psi_\lambda (z,y), \oPsi_\lambda (z,y)\]
%exists and are analytic.

Since $[e_n,z] =z_n$, the condition $(1+[e_n,z])/2 \not\in [1,\infty)$ is equivalent to 
$z_n\in \C{} \setminus [1,\infty) $. In particular, we have
$y\mapsto \Psi_\lambda (e_n,y)$ is analytic on $\{y\in \dS \: y_n< 1\}$. The same holds for
the functions $y\mapsto \Psi_\lambda (y,e_n)$, $\oPsi_\lambda (e_n,y)$ and $\oPsi_\lambda (y,e_n)$.

A vector $v \in \R{1,n}$ is {\it time-like} if $[v,v]<0$, {\it space-like} if $[v,v]>0$ and   {\it light-like}   if $[v,v]=0$.
For $x \in \dS$, we denote the future, past  and light cones of $x$ by $\Gamma^{+}(x)$, $\Gamma^-(x)$, 
and $\Gamma^0(x)$, respectively, where
\[
\Gamma^{\pm}(x)  := \{y\in \dS\: [y-x,y-x] <0, \pm y_0 >0\},\]
\[\Gamma (x) := \{y\in \dS\: [y-x,y-x]<0\} = \Gamma^+(x)\cup\Gamma^-(x)
\]
and 
\[\Gamma^0(x) := \{y \in \dS: [y-x,y-x]=0\}.\]

\begin{lemma}\label{pro:Gamma}
    For $y\in \dS$, $y_n \geq 1 \Leftrightarrow y \in \cl{\Gamma(e_n)}$. \end{lemma}
\begin{proof} We have
    \begin{align*}
        [y-e_n, y-e_n] &= [y,y] -2[y,e_n] +[e_n,e_n]\\
        &= 1 -2y_n +1 = 2(1-y_n)
    \end{align*}
and $2(1-y_n)\le 0 $ if and only if $y_n\ge 1$. 
\end{proof}

\begin{lemma}\label{lem:gGamma} Let $x\in \dS$ and $g\in G$. Then $\Gamma (gx) = g\Gamma (x)$.
\end{lemma}
\begin{proof} We have $[y - gx , y -gx] = [g^{-1} y-x,g^{-1}y-x]$ because of the invariance of
$[\cdot , \cdot]$. This implies that $y \in \Gamma (g x)$ if and only if $g^{-1}y\in \Gamma (x)$, which is
equivalent to the condition $y \in g\Gamma (x)$.
\end{proof}
\begin{theorem}\label{thm:4.5}  Let $x\in \dS$. Then  $\ltfrac{1+[y,x]}{2} \in [1,\infty)$ if and only if $y \in \cl{\Gamma(x)}$. 
\end{theorem}
\begin{proof} Write $x =ge_n$. Then $[y,x] = [y,ge_n] = [g^{-1}y,e_n]$.
The claim now follows from Lemma \ref{pro:Gamma} and Lemma ~\ref{lem:gGamma}. 
\end{proof}
For $f,g$ positive functions, we write $f \asymp g$ if there exist $M_1,M_2>0$ such that $M_1 g \leq f \leq M_2 g$. The following result follows immediately from Theorem~\ref{thm:hgf} on the asymptotic behavior of the function $\Fl(z)$ and the preceding theorem.
\begin{theorem}\label{cor:asymp}
    Let $x = ge_n$. In a small open neighborhood of the set $\{y \in \dS: [x,y] = 1\}$, we have the following asymptotic behavior near $t = \pi/2$: 
    \[\Psi_\lambda(gz_t,y) \asymp -\ltfrac{1}{\Gamma(\frac{1}{2} + \lambda)\Gamma(\frac{1}{2} - \lambda)} \log \left(\frac{1-[gz_t,y]}{2}\right),\]
     for $n=2$ and 
    \[\Psi_\lambda(gz_t,y) \asymp \ltfrac{\Gamma(n/2)\Gamma((n-2)/2)}{\Gamma(\rho +\lambda)\Gamma(\rho-\lambda)} \Bigg(\ltfrac{1-[gz_t,y]}{2}\Bigg)^{\frac{2-n}{2}},\]
    for $n \geq 3$. We obtain the same asymptotic behavior for $\oPsi_\lambda$ by taking complex conjugation.
\end{theorem}

\subsection{The distributions $\Psi_{\lambda,x}$ and $\oPsi_{\lambda, x}$}
For $\lambda \in \C{}$, let $a= \rho + \lambda$, $b =\rho-\lambda $ and $c= n/2$. We will simplify the notation by writing  $\Fl(z)$
instead of ${}_2F_1(a,b;c;z)$.  Formally, for  $x = g e_n \in \dS$, we define 
\begin{equation}\label{eq:psidist}
    \Psi_{\lambda,x} := \lim_{t \nearrow \pi/2}\Psi_\lambda(gz_t,y),\quad
    \oPsi_{\lambda,x } := \lim_{t \nearrow \pi/2}\Psi_\lambda(g\overline{z_t},y),
\end{equation}
taken in the sense of (\ref{eq:theta limit}). 

The following lemma shows that the limits defined above in (\ref{eq:psidist}) exist in the space of distributions for $x = g\cdot e_n$. However, it will be proved in \cref{sec:proof} that if $x = g_1 e_n = g_2 e_n$, then $\Psi_{\lambda,x} = \lim_{t \nearrow \pi/2}\Psi_\lambda(g_1z_t, \cdot) = \lim_{t \nearrow \pi/2}\Psi_\lambda(g_2z_t, \cdot)$ in the space of distributions. This forms part of the proof of Theorem~\ref{thm:1}. Thus, the limits taken in the sense of (\ref{eq:theta limit}) are well-defined.
\begin{lemma}\label{lemma:psi} Let $\lambda \in \C{} \setminus \{(\rho + \mathbb{N}_0) \cup (-\rho - \mathbb{N}_0)\}$ and 
$g \in \mr{SO}_{1,n}(\R{})_e$. Then the following holds:
\begin{enumerate}[leftmargin=*,label={\arabic*)}]
\item For $n \geq 2$, the limits
\[\Psi_{\lambda,x} = \lim_{t \nearrow\pi/2}\Psi_\lambda(g.z_t,\,\cdot \, )\quad \text{and}\quad \overline{\Psi}_{\lambda,x}
=
\lim_{t \nearrow\pi/2}\oPsi_\lambda(g.\bar{z}_t,\,\cdot \, ),\quad x= g\cdot e_n,\] exist in the space of distributions on $\dS$. \\

\item The distributions $\Psi_{\lambda,x}$ and $\oPsi_{\lambda,x}$  satisfy the  differential equation 
\[(\Delta - (\rho^2-\lambda^2))\Psi_{\lambda,x} =   (\Delta - (\rho^2-\lambda^2))\overline{\Psi}_{\lambda,x}=0.\]
\item  Moreover, $\Psi_{\lambda,x}$ and $\overline{\Psi}_{\lambda,x}$ can be represented by analytic functions in the following regions:
\begin{align*}
	\Psi_{\lambda,x}(y) = \begin{cases}
		{}_2F_1\Big(\frac{1+[x,y]}{2}\Big) &\text{if $y \notin \overline{\Gamma(x)}$},\\
		{}_2F_1\Big(\frac{1+[x,y]}{2} - i0\Big)  &\text{if $y \in \Gamma^+(x)$}, \\
		{}_2F_1\Big(\frac{1+[x,y]}{2} + i0\Big)  &\text{if $y \in \Gamma^-(x)$};
	\end{cases}\\[2mm]
	\overline{\Psi}_{\lambda,x}(y) = \begin{cases}
		{}_2F_1\Big(\frac{1+[x,y]}{2}\Big) &\text{if $y \notin \overline{\Gamma(x)}$},\\
		{}_2F_1\Big(\frac{1+[x,y]}{2} + i0\Big)  &\text{if $y \in \Gamma^+(x)$}, \\
		{}_2F_1\Big(\frac{1+[x,y]}{2} - i0\Big)  &\text{if $y \in \Gamma^-(x)$}.
	\end{cases}  
\end{align*} 
\end{enumerate} 
  
\end{lemma}

\begin{proof} For  $g \in G$, $x=g  e_n$ and $\varphi \in C_c^\infty(\dS)$, we have
\begin{align*}
    \Psi_{\lambda,x}(\varphi) &=
    \lim_{t \nearrow \pi/2}  \int_{\dS}\Psi_\lambda(gz_t,y) \overline{\varphi(y)}dy\\ 
   & = \lim_{t \nearrow \pi/2} \int_{\dS}\Psi_\lambda(z_t,g^{-1}y) \overline{\varphi(y)}dy.
   \end{align*}
  Since $dy$ is a left $G$-invariant measure, we can change variables and get
  \[\Psi_{\lambda,x}(\varphi)
    = \lim_{t \nearrow \pi/2}\int_{\dS}\Psi_\lambda(z_t,y) \overline{\varphi(gy)}dy
     =\Psi_{\lambda, e_n} (\varphi \circ\ell_g).\]
Similarly, we have $\oPsi_{\lambda,x}(\varphi) = \oPsi_{\lambda, e_n}(\varphi \circ\ell_g)$. 
 Thus, it is enough to prove the theorem for $\Psi_{\lambda,e_n} $. 
 
 Let $V_{x} = \kappa_{x}U_x$ be the local coordinate chart around $x$ as in Lemma \ref{lemma: exp_x}. Using standard
partition of unity  arguments together with  
\[ \lim_{t \nearrow \pi/2} \int_{V_x} \Psi_{\lambda}(z_t,y)\overline{\varphi(y)} dy = \lim_{t \nearrow \pi/2}\int_{V_{e_n}} \Psi_{\lambda}(z_t,g^{-1}y)\overline{\varphi(g^{-1}y)} dy, \]
where $V_{e_n} = g^{-1}V_x$,  we can assume that $\mathrm{supp}( \varphi) \subset V_{e_n}$.  
 With the same argument as above, it is enough to prove the existence of
  $\lim_{t \nearrow \pi/2}\kappa_{e_n}^*\oPsi_\lambda(\overline{z_t},v)$ in the space of distributions on $U_{e_n}$.

The distribution  \[\kappa_{e_n}^*\Psi_\lambda(z_t, v) = \hg{1 - iS([v,v])v_0\cos t + C([v,v])\sin t}\]
is real analytic as a function of $v \in U_{e_n}$ for $t < \pi/2$. Moreover, from Theorem~\ref{thm:hgf} we have the following pointwise limit: 
\[\kappa_{e_n}^*\Psi_{\lambda, e_n}(v)=\lim_{t \nearrow \pi/2} \kappa_{e_n}^*\Psi_\lambda(z_t, v) = \begin{cases}
    \hg{1+C[v,v]},  &  [v,v] >0,\\
    \hgf{\frac{1+ C[v,v]}{2} - i0},  &  [v,v] < 0, v_0 >0,\\
    \hgf{\frac{1+ C[v,v]}{2} + i0},  &  [v,v] < 0, v_0 <0.
    \end{cases}\]

By Theorem~\ref{thm:hgf}, the limit is uniform on compact sets in the disjoint regions $\{v:[v,v] >0\}$, $\{v:[v,v] < 0, v_0 >0\}$  and $\{v:[v,v] < 0, v_0 <0\}$. Also, $\kappa_{e_n}^*\Psi_{\lambda, e_n}(v)$ is real analytic in these disjoint open  regions. Thus, for  $n \geq 2$, if $\varphi$ is such that $\mathrm{supp}(\varphi) \cap \{v:[v,v] = 0\} = \emptyset $, then we have 
$$ \underset{t \nearrow \pi/2}{\mathrm{lim}}\;\int_{U_{e_n}} \kappa_{e_n}^*\Psi_{\lambda}(z_t, v) \overline{\varphi(v)} dv \longrightarrow \int_{U_{e_n}} \kappa_{e_n}^*\Psi_{\lambda, e_n}(v) \overline{\varphi(v)} dv, $$
 where $dv = \kappa^* d\mu_{\dS}(y)$ denotes the pullback to $U_{e_n}$ of the restriction of the $G$-invariant measure $d\mu_{\dS}$ to $V_{e_n}$. 
 
In the following, let $W_{e_n } \subset U_{e_n}$ be a sufficiently small open neighborhood of the set $\{v \in \rT_{e_n}(\dS):[v,v] =0\}$ such that the asymptotic behavior described in Corollary~\ref{cor:asymp} holds.
\smallskip

\noindent
{\bf Case: $n=2$.}
From (\ref{eq:3}), Corollary~\ref{cor:asymp}, Appendix~\ref{appendixa} and  \cite[Chap. II, Sec. 2.4, Ex. 4]{GS64}, we obtain for $v \in W_{e_n}$: 
\[
    |\kappa_{e_n}^*\Psi_\lambda(z_t,v)| \asymp \left|\ltfrac{1}{\Gamma(\frac{1}{2} + \lambda)\Gamma(\frac{1}{2} - \lambda)} \ln\left(\frac{1-[z_t,S([v,v])v+C[v,v]e_n]}{2}\right)\right|,\]
    which implies that
   \begin{equation}\label{eq:estim}
   |\kappa_{e_n}^*\Psi_\lambda(z_t,v)| \leq d \left| 1+\ln\left|{1-C[v,v]}\right|\right|\end{equation}
   for some constant $d>0$.
 Since the logarithm is locally integrable, it follows that for $\Supp \varphi \subset W_{e_n}$, 
\[\int |\kappa_{e_n}^* \Psi_\lambda (z_t,v)\kappa_{e_n}^*\varphi (v)| dv \le d\int \left| \left(1+\ln\left|{1-C[v,v]}\right|\right)
\kappa_{e_n}^*\varphi (v)\right|
dv.\]
As the right hand side is independent of $t$, it follows by LDCT that the limit  $t \nearrow \pi/2$ also exists and defines a 
distribution.
\smallskip

\noindent
{\bf Case: $n \geq 3$.}
 We have $\mathrm{Re}(c-a-b) = \ltfrac{2-n}{2} < 0$ for $n \geq 3$. From Corollary~\ref{cor:asymp} and  (\ref{4}), we have the following estimates in the small neighborhood $W_{e_n}$ of the set $\{v:[v,v]=0\}$:
\begin{align*}
   \kappa_{e_n}^*\Psi_\lambda(z_t,v) &\asymp \ltfrac{\Gamma(n/2)\Gamma((n-2)/2)}{\Gamma(\rho +\lambda)\Gamma(\rho-\lambda)} \left(\ltfrac{1-[z_t,S([v,v])v+C([v,v])e_n]}{2}\right)^{\frac{2-n}{2}},\\ 
  | \kappa_{e_n}^*\Psi_\lambda(z_t,v)|& \leq d  \left|{1-C([v,v])}\right|^{\frac{2-n}{2}} \leq d |[v,v]|^{\frac{2-n}{2}}
\end{align*}
for some constant $d>0$.
Since the right hand side of the above inequality defines a distribution (see \cite[Chap. III, Sec. 2.8,2.9]{GS64}),  $\lim_{t \nearrow \pi/2}\kappa_{e_n}^*\Psi_\lambda(z_t,\cdot)$ exists in the space of distributions. The arguments from the beginning of the proof then yields the result for $\Psi_{\lambda, x}$ and $\overline{\Psi}_{\lambda, x}$.

 We now claim that $(\Delta - (\rho^2-\lambda^2))\Psi_{\lambda, e_n} = 0$. Using the fact that differentiation is a continuous linear map on the space of distributions, we obtain  $$\lim_{t \nearrow \pi/2}(\Delta - (\rho^2-\lambda^2))\Psi_\lambda(z_t, y) = (\Delta - (\rho^2-\lambda^2))\Psi_{\lambda, e_n}. $$
Now let $a = \rho +\lambda$, $b=\rho-\lambda$, $c = n/2$ and $w_t = \frac{1+[z_t,y]}{2} $. Using local coordinates and the properties of the hypergeometric function obtained by combining \cite[ (9.2.2)]{LS} and \cite[(9.2.16)]{LS}, we arrive at 
\begin{small}
    \begin{align*}
     (\Delta - (\rho^2-\lambda^2))\Psi_\lambda(z_t,y) & = \ltfrac{ab}{c}[w_t(1-w_t)\frac{(a+1)(b+1)}{(c+1)}{}_2F_1(a+2,b+2,c+2,w_t) \\
     &+(\frac{n}{2}-nw_t){}_2F_1(a+1,b+1,c+1,w_t)- c \,{}_2F_1(a,b,c,w_t)] \\
     &=0.
 \end{align*}
\end{small}

Therefore, as distributions, we have $(\Delta -(\rho^2-\lambda^2))\Psi_{\lambda, e_n} = 0$. Following the same steps, we obtain $(\Delta - (\rho^2-\lambda^2))\overline{\Psi}_{\lambda, e_n} = 0$. As $(\Delta - (\rho^2-\lambda^2))$ is a $G$-invariant operator, we also have \[(\Delta - (\rho^2-\lambda^2))\Psi_{\lambda,x} =  (\Delta - (\rho^2-\lambda^2))\overline{\Psi}_{\lambda,x}=0.\]
Lastly,  we observe from Theorem~\ref{thm:hgf} that the pointwise limit can be obtained as follows:
\begin{align*}
    \Psi_{\lambda, e_n}(y) &= \lim_{t \nearrow \pi/2} \hg{1+[z_t,y]}\\
    &=  \lim_{t\nearrow \pi/2}  \hg{1+\sin (t) y_n - i\cos (t)y_0} \\[2mm]
    &= \begin{cases}
        {}_2F_1\Big(\frac{1+y_n}{2}\Big) &\text{if $y_n < 1$},\\
        {}_2F_1\Big(\frac{1+y_n}{2} - i0\Big)  &\text{if $y_n > 1,y_0 > 0$}, \\
     {}_2F_1\Big(\frac{1+y_n}{2} + i0\Big)  &\text{if $y_n>1, y_0 < 0$},
    \end{cases}
\end{align*}
where $\hgf{x \pm i0}$ has been calculated in Appendix \ref{sec:hgf} for $x>1$. For $\Psi_{\lambda,x}(y)$ we use the relation $\Psi_{\lambda,x}(y) = \Psi_{\lambda,e_n}(g^{-1}y)$ for $x = ge_n$.  For $\overline{\Psi}_{\lambda,x}$, one gets the complex-conjugate expressions by Theorem~\ref{thm:kernel}. 
\end{proof}

\begin{Remark}(1) If $\lambda \in  (\rho + \mathbb{N}_0) \cup (-\rho - \mathbb{N}_0)$, then the kernels $\Psi_\lambda$ and $\oPsi_\lambda$ are polynomials, since either the first argument $a$ or the second argument $b$ of the hypergeometric function is a non-positive integer. Thus, the distributional limits $\Psi_{\lambda,x}$ and $\overline{\Psi}_{\lambda,x}$ are represented by polynomials in $[x,y]$ and $\Psi_{\lambda, x}= \oPsi_{\lambda,x}$.\\
	(2) If $\lambda \in i\R{} \cup (-\rho,\rho)$, then $\Psi_\lambda$ and $\overline{\Psi}_\lambda$ are the positive definite kernels discussed in Section~\ref{sec:ker}.
\end{Remark}

The distribution $\Psi_{\lambda,x} - \overline{\Psi}_{\lambda,x}$, which is supported only on the future and past light cones of $x \in \dS$, is studied in quantum field theory (see \cite{BJM23,BV96,BV97, BM96, BM04}). 

\begin{theorem}
Write 
\[F(\lambda, n; x,y)= {}_2F_1\left(\frac{1}{2}-\lambda, \frac{1}{2} + \lambda , \frac{4-n}{2}; \frac{1-[x,y]}{2}\right) \]
and 
\[ \Phi_{\lambda,x} = \Psi_{\lambda,x} - \overline{\Psi}_{\lambda,x}. \]
Then the following holds: 
\begin{enumerate}[leftmargin=*,label={\arabic*)}]
	\item For $n$ odd and $c_1 = (-1)^{\frac{n-1}{2}}\frac{i\Gamma(n/2)\Gamma((n-2)/2)}{\Gamma(\rho + \lambda)\Gamma(\rho-\lambda)}$:
	\[
	\Phi_{\lambda,x}(y)= 
	c_1\begin{cases}
		0 , & y \notin \cl {\Gamma(x)},\\
		\left(\frac{[x,y]-1}{2}\right)^{\frac{2-n}{2}} F(\lambda,n;x,y) ,
		& y \in \Gamma^+(x),\\
		-\left(\frac{[x,y]-1}{2}\right)^{\frac{2-n}{2}} F(\lambda, n;x,y), &y \in \Gamma^-(x).
	\end{cases}
	\]
	\item  For even $n$ and $c_2 = (-1)^{\frac{n-2}{2}} \frac{i\pi\Gamma(n/2)}{\Gamma(1/2 + \lambda)\Gamma(1/2-\lambda)}$ we have:
	\[
	\Phi_{\lambda,x}(y) =
	c_2\begin{cases}
		0 , & y \notin \cl {\Gamma(x)},\\
		\sum_{k=0}^{\infty}\frac{(\rho+\lambda)_k(\rho - \lambda)_k}{k!(n/2 - 1 +k)!}\left(\frac{[x,y]-1}{2}\right)^k , &  y \in \Gamma^+(x), \\
		-\sum_{k=0}^{\infty}\frac{(\rho+\lambda)_k(\rho - \lambda)_k}{k!(n/2 - 1 +k)!}\left(\frac{[x,y]-1}{2}\right)^k , & y \in \Gamma^-(x).
	\end{cases}   
	\] 
	\end{enumerate} 

\noindent
{\rm  (Huygens' Principle):} Assume that $n\geq 4$.  If
\[\lambda \in \Bigg\{\frac{1+2k}{2}: k=\lfloor -n/2+1\rfloor,...,\lfloor n/2-2\rfloor\Bigg\} \cap (-\rho ,\rho),\] 
then $\Phi_{\lambda,x}$ is supported on the boundary $\partial \overline{\Gamma(x)}$ of the light cone at $x$.

\end{theorem}

Note that if $n$ is odd, then the constant $c_1 \neq 0$ for any $\lambda \in i\R{} \cup (-\rho,\rho)$. However, if $n \geq 4$ and $n$ is even, then for the values of $\lambda$ in the set described above in the theorem, we have that the distribution $\Phi_{\lambda,x}  $ is zero on the open regions $\Gamma^+(x)$, $\Gamma^-(x)$ and $(\overline{\Gamma(x)})^c$. This phenomenon was also observed in \cite{FNO24}.

%---------------------------------------------------------------------------------------------------------------------------------------------------------------------------------------------------------------------------------------------------------------------------------------------------------------------------------------------------
%\vspace{5mm}
\section{Wavefront Sets}\label{sec: wf}

The wavefront set of a distribution was introduced by L. H\"ormander in 1970. It gives information about the singular support of a distribution and the directions in which the distribution is not smooth or analytic. We have seen in the previous section that the  distributions $\Psi_{\lambda,x}$ and $\overline{\Psi}_{\lambda,x}$ can be written as analytic functions everywhere on the de Sitter space except at the boundary of the light cone of $x$. We will study the singularities of these distributions in terms of wavefront sets. Let us now recall the definition of wavefront set. Our main reference is \cite[Chap. 8]{H03}.

Let $X \subset \R{1,n}$ be an open subset and suppose $\Theta \in \mathcal{E}'(X)$ is an anti-linear distribution with compact support. For $\xi\in \R{1,n}\setminus \{0\}$, we can define the Fourier transform of $\Theta$ as follows:
\begin{equation}
    \widehat{\Theta}(\xi) = \Theta(e^{-2\pi i [x,\xi]}),
\end{equation}
where $[x,\xi] = -x_0\xi_0 + x_1\xi_1+...+x_{n}\xi_{n}$. We use $[\cdot,\cdot]$ instead of the usual scalar product to ensure $G',G$-covariance.

\begin{definition}\label{def:wfsmooth}
    Let $\Theta \in \mathcal{D}^*(X)$. We say that $(x_0,\xi_0) \in \mathrm{T}^*(X)\setminus (X \times \{0\})$ is a regular directed point if there exist an open neighbourhood $U$ of $x_0$, a conical neighbourhood $V$ of $\xi_0$ and $\varphi \in C_c^{\infty}(U)$ with $\varphi(x_0) \neq 0$ such that for all $N \in \mathbb{N}$ and $\xi \in V$:
\begin{equation}\label{ft}
    |\widehat{\varphi \Theta}(\xi)| \leq C_{N,\varphi}(1+|\xi|)^{-N}.
\end{equation}
 The wavefront set $WF(\Theta) \subseteq \mathrm{T}^*(X) \setminus (X \times \{0\})$ is the complement of the regular directed set.
\end{definition}

\begin{Remark}Let $\Theta \in \mathcal{D}^*(X)$. Then $\overline{\Theta}$ is a linear distribution. Let $\mathscr{F}$ be the Euclidean Fourier transform. We will denote the wavefront sets obtained using the Euclidean inner product in the definition of the Fourier transform by ${WF}_{\mr{Euc}}$.  Let $J = \mr{diag}(-1,1,...,1)$ be an $(n+1) \times (n+1)$ matrix. For a smooth function $\varphi$ with compact support, we have 
    \begin{align*}
        \widehat{\varphi \Theta}(\xi) &= \Theta_x (\varphi(x)e^{-2\pi i [x,\xi]}) = \Theta _x(\varphi(x)e^{-2\pi i \langle x,J\xi\rangle})\\
        & = \overline{\mathscr{F} (\overline{\varphi \Theta})(-J\xi)},
    \end{align*}
    where $\langle \cdot, \cdot \rangle$ denotes the Euclidean inner product. It then follows that
\[(x_0,\xi_0) \notin WF(\Theta) \Longleftrightarrow(x_0,-J\xi_0) \notin WF_{\mr{Euc}}(\overline{\Theta}).\]    
Thus, all the results proved by H\"ormander in \cite{H03} using the Euclidean inner product for linear distributions also hold for antilinear distributions when the Lorentzian bilinear form is used in the definition of the Fourier transform. Hence, we will state the results from \cite[Chap. 8]{H03} for the Lorentzian bilinear form.
\end{Remark}

\begin{definition}
    Let $\Theta \in \mathcal{D}^*(X)$. The {\it singular support} of $\Theta$ is the set of all points $x$ for which there exists no neighbourhood of $x$ on which the restriction of  $\Theta$ is a $C^{\infty}$ function.
\end{definition}

\begin{lemma}
    If $\Theta \in \mathcal{D}^*(X)$, then the projection of $WF(\Theta)$ onto $X$ is the singular support of $\Theta$.
\end{lemma}

\begin{Remark}
    The wavefront set $WF(\Theta)$ is conic; that is, if $(x,\xi) \in WF(\Theta)$, then for every $s > 0$, we also have $(x,s \xi) \in WF(\Theta)$.
\end{Remark}

We will now define the analytic wavefront set and follow the definition in \cite[Def. 8.4.3]{H63}. As there are no compactly supported analytic functions, we cannot use Definition~\ref{def:wfsmooth} to define the analytic wavefront set.

\begin{definition}
    Let $X$ be an open subset of $\R{1,n}$ and $\Theta \in \mathcal{D}^*(X)$. We define $WF_A(\Theta)$ to be the complement in $X \times (\R{1,n}\setminus \{0\})$ of the set of points $(x_0,\xi_0)$ such that there exists an open neighbourhood $U \subset X$ of $x_0$, a conic neighbourhood $\Gamma$ of $\xi_0$ and a bounded sequence of $\Theta_N \in \mathcal{E}'(X)$ that are equal to $\Theta$ in $U$ and satisfy 
$$||\xi|^N\widehat{\Theta_N}(\xi)| \leq C(C(N+1))^N ,\qquad N=1,2,...,$$
for all $\xi \in \Gamma$ and for some $C>0$.
\end{definition}

\begin{example}
    Let $u=\delta_0$ in $\R{1,n}$. We will use \cite[Lem. 8.4.4]{H63} to determine the analytic wavefront set. Observe that $WF_A(\delta_0) \subset \{0\}\times (\R{1,n}\setminus \{0\})$. Let $\chi_N$ be a sequence of functions as in \cite[Lem. 8.4.4]{H63}. For $\xi \neq 0$, we have
    \begin{align*}
        ||\xi|^N\widehat{\chi_N\delta_0}(\xi)| = ||\xi|^N\chi_N(0) |= |\xi|^N,
    \end{align*}
    which is not bounded as in Definition 5.6. Therefore, $WF_A(\delta_0) = \{(0,\xi): \xi \in \R{1,n}\setminus \{0\} \}. $ \qed
\end{example}
The following lemma describes the relation between the smooth wavefront set and the analytic wavefront set of a distribution.
\begin{lemma}
    Let $\Theta \in \mathcal{D}^*(X)$. Then  $WF(\Theta) \subset WF_A(\Theta)$. 
\end{lemma}

    \begin{proof}
        Suppose $(x_0,\xi_0) \notin WF_A(\Theta)$. Then there exist a  neighborhood $U \ni x_0$, a closed cone $\Gamma \ni \xi_0$ and a sequence of functions $\chi_N \in C_c^\infty(U)$  as in \cite[Lem. 8.4.4]{H63} such that  
        $$||\xi|^N\widehat{\chi_N\Theta}(\xi)| \leq C(C(N+1))^{N}, \quad \xi \in \Gamma.$$ 
        Thus, it follows from the definition that $(x_0, \xi_0) \notin WF(\Theta)$. 
    \end{proof}

\par

The following result, which is an analogue of \cite[Thm. 8.2.4]{H63}, is a rephrased version of \cite[Thm. 8.5.1]{H63}, stated in a form parallel to that of \cite[Thm. 8.2.4]{H63}.

\begin{theorem} \label{cw}
    Let $X \subset \mathbb{R}^n$ and $Y \subset \mathbb{R}^m$ be open sets and let $\iota : X \to Y$ be a real-analytic map. Denote the normal set of the map by
    $$N_{\iota} \:= \{(\iota(x),\xi) \in Y\times \R{n} : {}^td\iota_x(\xi) =0 \}.$$ 
     If $\Theta \in \mathcal{D}^*(Y)$  is such that $N_\iota \cap WF_A(\Theta) = \emptyset$, then the pullback $\iota^*\Theta$ is well defined and 
    $$WF_A(\iota^*\Theta) \subset \iota^*(WF_A(\Theta)),$$ 
    where  $$\iota^*(WF_A(\Theta))= \{(x,{}^td\iota_x(\xi)): (\iota(x),\xi) \in WF_A(\Theta)\}.$$
\end{theorem}

Theorem~\ref{cw} shows that the definition of the analytic wavefront set is invariant under coordinate changes and thus allows us to define the analytic wavefront set on a real analytic manifold $X$. 

\begin{definition}
    Let $X$ be a real analytic manifold and let $(U_\kappa,\kappa)$ be analytic local coordinates on X. We define $WF_A(\Theta) \subset \mathrm{T}^*(X)\setminus (X \times \{0\})$ to be the set \[\bigcup_{\kappa^*}WF_A((\kappa^{-1})^*\Theta) :=\bigcup_{\kappa} \{(x, {}^td\kappa^{-1}_x(\eta)): (\kappa^{-1}(x),\eta) \in WF_A((\kappa^{-1})^*\Theta)\},\] where $(\kappa^{-1})^*\Theta (\varphi)= \Theta(\varphi\circ \kappa^{-1})$ for any $\varphi \in C_c^\infty (U_\kappa)$.
\end{definition}

The next theorem describes the analytic wavefront set of distributions which are boundary values of analytic functions. Let $\Gamma$ be an open convex cone in $\R{n+1}$. Then the {\it dual cone} $\Gamma^{\circ}$ is defined as 
$$\Gamma^{\circ} = \{ \eta \in \R{n+1} : (\forall \xi \in \Gamma) \; [\eta,\xi]  \leq 0\}.$$
\begin{theorem}\label{bd}
Let $X \subset \R{n+1}$ be an open set and $\Gamma$ an open convex cone in $\R{n+1}$. For some $\gamma > 0 $, let
$$Z = \{ z \in \C{n+1} : {\rm Re}\, z \in X, \, {\rm Im}\, z \in \Gamma, \, |{\rm Im}\, z| < \gamma\}.$$
If $\Theta$ is an analytic function in $Z$ such that 
$$|\Theta(z)| \leq C |{\rm Im} \,z|^{-N}$$
for some positive integer $N$ and some constant $C >0$, then  ${\lim}_{y\searrow 0}\Theta(.+iy) = \Theta_0$ exists in the sense of distributions and is of order $N$.  We also have 
$$ WF_A (\Theta_0) \subset X \times (\Gamma^{\circ} \setminus 0).$$
\end{theorem}
\begin{proof}
    See Theorem 3.1.15 and Theorem 8.4.8 in \cite{H63}.
\end{proof}
\begin{comment}

\begin{example}
    Consider the distribution on $\R{}$, $\Theta = (x+i0)^{\frac{2-n}{2}}$. It is the limit of the analytic function $(x+iy)^{\frac{2-n}{2}}$ for $x \in \R{}$ and $y \in \Gamma = \R{}_+$. Then its dual cone is $\Gamma^{\circ} = \R{}_{\geq 0 }$. By \cref{bd}, $WF_A(\Theta) \subset \R{} \times \R{}_+$. It is obvious that the distribution has singularities only at $x =0$. Therefore,
    $$WF_A((x+i0)^{\frac{2-n}{2}}) = \{(0, \tau): \tau > 0\}.$$
    Similarily,
    \begin{equation*}\pushQED{\qed}
        WF_A((x-i0)^{\frac{2-n}{2}}) = \{(0, \tau): \tau < 0\}. \qedhere
    \end{equation*}
    
\end{example}
\end{comment}
\begin{example}
    Let $\Theta = \ln({x+i0})$, which is a boundary value of the holomorphic function $\ln({x +iy})$ for $y >0$ (see \cite[p. 26, Ex. 6]{GS64}). Since the logarithm grows slower than any negative power of $|y|$, the limit as $y \searrow 0$ is a distribution on $\R{}$. It follows from Theorem~\ref{bd} that 
    $$WF_A(\ln({x+i0})) = \{(0,\xi):\xi>0\}. $$
    Likewise, we have
    \begin{equation*}\pushQED{\qed}
        WF_A(\ln({x-i0})) = \{(0,\xi): \xi<0\}. \qedhere
    \end{equation*}
    
\end{example}

\begin{example}\label{eg:hgf}
    Let $\Theta = {}_2F_1(x + i0) := \lim_{y \searrow 0}{}_2F_1(x +iy)$. We have proved in Appendix~\ref{sec:hgf} that it is a distribution which can be represented as a real analytic function everywhere except at $x=1$. As a result of Theorem~\ref{bd}, the analytic wavefront set is $$WF_A({}_2F_1(x + i0)) = \{(1, \xi):\xi > 0\},$$ 
    and it also follows that 
    \begin{equation*}\pushQED{\qed}
       WF_A({}_2F_1(x - i0))= \{(1, \xi):\xi < 0\}.\qedhere  
    \end{equation*}
  \end{example}

 Let $P(x,D) = \underset{|\alpha| \leq m}{\Sigma} a_\alpha(x) D^\alpha$ be a differential operator on $X$ with analytic coefficients. Then $$WF_A(P(x,D)\Theta) \subset WF_A(\Theta).$$
The following is a restatement of \cite[Thm. 8.6.1]{H63}, which is a converse to the above theorem.

\begin{theorem}\label{thm:char}
    If P(x,D) is a differential operator of order m with real analytic coefficients in X, then
    $$WF_A(\Theta) \subset WF_A(P\Theta) \cup \mr{Char} (P),$$
    where the characteristic set of P is defined by
    $$\mr{Char} (P) = \{ (x,\xi) \in T ^*(X)\setminus 0 : P_m(x,\xi) := \underset{|\alpha| = m}{\Sigma}a_\alpha \xi^\alpha=0\}.$$
\end{theorem}

Let $P(x,D)$ be a differential operator of order $m$ on a manifold $X$ with real analytic coefficients  and principal symbol $P_m = \underset{|\alpha| = m}{\Sigma}a_\alpha \xi^\alpha$ in local coordinates. We say that the curve $(x(t),\xi(t))$ in $\rT^*(\dS)$ is a bicharacteristic strip if $P_m(x(t), \xi(t))=0$ for all $t$, with initial data $(x_0,\xi_0) \in \mr{Char} (P_m)$ and it satisfies the Hamiltonian equations given by
\[ \frac{dx}{dt} = \frac{\partial P_m(x,\xi)}{\partial \xi}, \qquad  \frac{d\xi}{dt} =- \frac{\partial P_m(x,\xi)}{\partial x}.\]

Let $S$ be a closed conic set in $\rT^*(X)$. We say that it is {\it invariant under the  Hamiltonian flow of $P_m$} if $S \subseteq {\rm Char} (P_m)$ and for  a bicharacteristic strip $(x(t),\xi(t))$ passing through $(x_0,\xi_0) \in S$, which satisfies the condition that $\partial P_m(x_0,\xi_0)/\partial \xi \neq 0$, the curve $(x(t),\xi(t))$ must lie in $S$ for all $t$. \\

The following result is \cite[Thm. 7.1]{H71}.

\begin{theorem}[Propagation of Singularities] \label{PS} 
    Let $P(x,D)$ be a differential operator with real analytic coefficients of order $m$ in $X$ and let $P_m$ be its real principal symbol. If $\Theta \in \mathcal{D}^*(X)$ and $P(x, D)\Theta = f$, it follows that $WF_A(\Theta) \setminus WF_A(f)$  is invariant under the Hamiltonian vector field of $P_m$ when $\partial P_m(x,\xi)/\partial \xi \neq 0$.
\end{theorem}

We say that a curve $(x(t), \xi(t))$ is a {\it null geodesic strip} if $[\dot{x}(t),\dot{x}(t)] = 0$ and $\xi(t)$ is the dual of $\dot{x}(t)$. The following proposition is proved in \cite[Prop. 2.8]{R92} for a curved space-time. In particular, the result holds for the de Sitter space. 

\begin{proposition}
   On $\dS$, the bicharacteristic strips for $\Delta - (\rho^2-\lambda^2)$ are exactly the null geodesic strips for $\dS$.
\end{proposition}

%\vspace{5mm}

Using analytic wavefront sets, we can prove that the distributions discussed in Section 4 cannot vanish on any non-empty open set $O$ of $\dS$. The following theorem is due to Strohmaier, Verch and Wollenberg (see \cite[Prop. 5.3]{SVW02}).

\begin{theorem}\label{thm:supp}
    Let X be a real analytic manifold and let $\Theta \in \mathcal{D}^*(X)$ be such that $WF_A(\Theta) \cap -WF_A(\Theta) = \emptyset$. If $\mathcal{O}$ is an open region in $X$ and 
    $\Theta|_\mathcal{O}  = 0$, then it follows that $\Theta=0$. Here, $-WF_A(\Theta)$ is defined as $\{(x,\xi) : (x, -\xi) \in WF_A(\Theta)\}$.
\end{theorem}
\begin{Remark}
	This theorem is not true for the smooth wavefront set. Consider the distribution \[\Theta(x) = \begin{cases}
		e^{-1/x} & \text{if $x > 0$},\\
		0 & \text{if $x\leq 0$}.
	\end{cases}\] Since $\Theta$ is a smooth function,  the wavefront set is
$WF(\Theta) = \emptyset$ and satisfies the condition that $WF(\Theta) \cap -WF(\Theta) = \emptyset$. Clearly, $\Theta$ is not the zero distribution; however, it is zero in the open region $(-\infty,0)$.
\end{Remark}

%---------------------------------------------------------------------------------------------------------------------------------------------------------------------------------------------------------------------------------------------------------------------------------------------------------------------------------------------------
\section{Wavefront set of spherical distributions on $\dS$} \label{sec:ws}

In this section we will restate Theorem~\ref{thm:1} and study its implications. The proof of Theorem~\ref{thm:1} is postponed to the next section. Recall that 
\[\rT_x(\dS) = \{v \in \R{n+1} : [x,v] = 0\}.\]
 We denote by  $\Psi_{\lambda,x}$ the distributional limit from Section \ref{sec:bv}.

\begin{manualtheoremx}{1} Let $\rho =(n-1)/2$ and  $G = \SO_{1,n}(\R{})_e$. Let $z_t = i\cos t e_0 + \sin t e_n \in \Xi$ and $\overline{z_t}\in \oXi$ for $t\in (-\pi/2, \pi/2)$. Then the following holds
	for $\lambda \in \C{} \setminus \{(\rho + \mathbb{N}_0) \cup (-\rho - \mathbb{N}_0)\}\colon$
\begin{enumerate}[leftmargin=*,label={\arabic*)}]
	\item  For $x\in \dS$,  let $g_1,g_2 \in G$ be such that $x = g_1e_n=g_2e_n$. Then
	\[\Psi_{\lambda,x} = \lim_{t \nearrow \pi/2}\Psi_\lambda(g_1z_t,\cdot)=\lim_{t \nearrow \pi/2}\Psi_\lambda(g_2z_t,\cdot),\]
	and 
	\[\oPsi_{\lambda,x} = \lim_{t \nearrow \pi/2}\oPsi_\lambda(g_1\oz_t,\cdot)=\lim_{t \nearrow \pi/2}\oPsi_\lambda(g_2\oz_t,\cdot)
	\]
	exist in the space of distributions. In particular, these limits are well defined for $x \in \dS$. 
	
	\item The analytic wavefront sets of the distributions $\Psi_{\lambda,x}$ and $\oPsi_{\lambda,x}$ are given by 
	\begin{align*}
		WF_A(\Psi_{\lambda,x}) &= \{(x,v)\: v_0 < 0,\, [x,v]= [v,v]=0 \}\\
		& \bigcup \{(x +v, (-v_0, \bv) ) \: v_0<0,\,  [x,v]= [v,v]=0\}\\
		&\bigcup \{(x+v, (v_0, -\bv))\: v_0 <0, \,  [x,v]=[v,v]=0\} 
	\end{align*}
	and
	\begin{align*}
		WF_A(\oPsi_{\lambda,x}) &=
		\{(x,v)\: v_0 > 0,\, [x,v]= [v,v]=0 \}\\
		& \bigcup \{(x +v, (v_0, -\bv) ) \: v_0>0,\,  [x,v]= [v,v]=0\}\\
		&\bigcup \{(x+v, (-v_0, \bv))\: v_0 >0, \,  [x,v]=[v,v]=0\}   . 
	\end{align*} 
\end{enumerate}

\end{manualtheoremx}

\begin{figure}[h]
\centering
\begin{subfigure}{0.5\textwidth}
	\centering
\begin{tikzpicture}[scale = 0.5, rotate = 90]
\draw[->, thick] (0,5)--(0,-5) node[below]{${\bf v}$};
 \draw[blue, ultra thick] (-4,-4)--(4,4);
 \draw[blue, ultra thick] (-4,4)--(4,-4);
 \filldraw[color=blue, fill = blue, fill opacity = 0.2,ultra thick] (4,0) ellipse (0.5 and 4) ;
  \filldraw[color=blue, fill = blue, fill opacity = 0.2,ultra thick] (-4,0) ellipse (0.5 and 4) ;
 \fill[ color=blue!40, opacity=0.2] (4,4)--(0,0)--(4,-4) arc (270: 90:0.5 and 4);
\fill[ color=blue!40, opacity=0.2] (-4,4)--(0,0)--(-4,-4) arc (270:90:0.5 and 4);
\fill[color=red, opacity=0.2] (1,1)--(0,0)--(1,-1) arc (270:90:0.2 and 1);
 \filldraw[color=red, fill = red,fill opacity = 0.2,ultra thick] (1,0) ellipse (0.2 and 1) ;
 \draw[->, thick] (-5,0)--(5,0) node[right]{$v_0$};
 \draw[->, red, thick] (2,2)--(3,1);
 \draw[->,red,thick] (2,-2)--(3,-1);
 \draw[->,red,thick] (-2,2)--(-1,3);
 \draw[->,red,thick] (-2,-2)--(-1,-3);
 \draw[red,thick](0,0)--(1,1);
 \draw[red,thick](0,0)--(1,-1);
\end{tikzpicture}
\caption{Analytic wavefront set of $\overline{\Psi}_{\lambda,x}$.}
\end{subfigure}

\par\bigskip 
\begin{subfigure}{0.5\textwidth}
    \centering
   \begin{tikzpicture}[scale = 0.5, rotate = 90]
\draw[->,thick] (0,5)--(0,-5) node[below]{${\bf v}$};
 \draw[blue, ultra thick] (-4,-4)--(4,4);
 \draw[blue, ultra thick] (-4,4)--(4,-4);
 %\filldraw[color=red!60, fill=red!5, very thick](-1,0) circle (1.5)
 %(2.5,0) ellipse (1.5 and 0.5);
 \filldraw[color=blue, fill = blue, fill opacity = 0.2,ultra thick] (4,0) ellipse (0.5 and 4) ;
  \filldraw[color=blue, fill = blue, fill opacity = 0.2,ultra thick] (-4,0) ellipse (0.5 and 4) ;
\fill[ color=blue!40, opacity=0.2] (4,4)--(0,0)--(4,-4) arc (270: 90:0.5 and 4);
\fill[ color=blue!40, opacity=0.2] (-4,4)--(0,0)--(-4,-4) arc (270:90:0.5 and 4);
\fill[color=red, opacity=0.2] (-1,1)--(0,0)--(-1,-1) arc (270:90:0.2 and 1);
 \filldraw[color=red, fill = red,fill opacity = 0.2,ultra thick] (-1,0) ellipse (0.2 and 1) ;
 \draw[->, thick] (-5,0)--(5,0) node[right]{$v_0$};
 \draw[->, red, thick] (2,2)--(1,3);
 \draw[->,red,thick] (2,-2)--(1,-3);
 \draw[->,red,thick] (-2,2)--(-3,1);
 \draw[->,red,thick] (-2,-2)--(-3,-1);
 \draw[red,thick](0,0)--(-1,1);
  \draw[red,thick](0,0)--(-1,-1);
\end{tikzpicture}
\caption{Analytic wavefront set of $\Psi_{\lambda,x}$ .}
\end{subfigure}
\caption{In this figure, the tangent space at a point $x$ has been identified with its cotangent space at $x$. The blue region is the light cone of 0, while the red arrows and the red region lie in the cotangent space at that point. The light cone, together with the red arrows and the red region, constitutes the analytic wavefront set.}
\label{fig:wavefront set}
\end{figure}
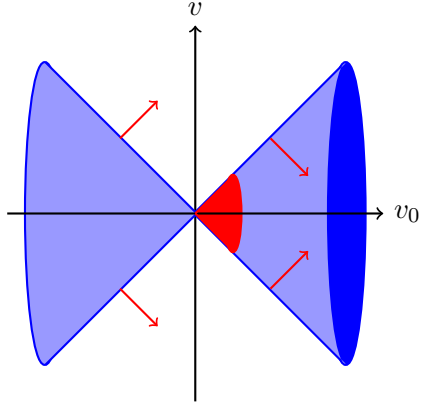
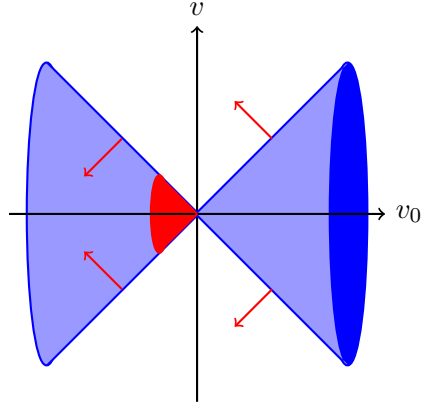

Observe that $[x+v,x+v] = [x,x] + 2[x,v] +[v,v] =1$, since $[x,x] =1$, $[x,v]=0$ and $[v,v]=0$. Therefore, $x + v \in \dS$.
Figure~\ref{fig:wavefront set} shows the analytic wavefront sets of $\Psi_{\lambda,x}$ and $\overline{\Psi}_{\lambda, x}$ viewed locally in $\rT_x(\dS)$.

\begin{corollary}\label{cor: pmwf} The wavefront sets of the distributions $\Psi_{\lambda,x} \pm \overline{\Psi}_{\lambda,x}$ are 
\begin{align*}
     WF_A(\Psi_{\lambda,x} \pm \overline{\Psi}_{\lambda,x}) &= WF_A(\Psi_{\lambda,x})\sqcup WF_A(\overline{\Psi}_{\lambda,x})\\
     &=\{(x,v): v_0 \neq 0\} \cup \{(x+v, (-v_0,\bv) ): v_0 \neq 0 \},
\end{align*}
where $v \in \rT_x(\dS)$ satisfies $[v,v]=0$.

\end{corollary}
\begin{proof} From Theorem~\ref{thm:1}, observe that the sets $WF_A(\Psi_{\lambda,x})$ and $WF_A(\overline{\Psi}_{\lambda,x})$ are mutually disjoint. Let us write $\Theta_{\pm} = \Psi_{\lambda,x} \pm \overline{\Psi}_{\lambda,x}$. 
   It is immediate from the definition of analytic wavefront sets that $WF_A(\Psi_{\lambda,x} \pm \overline{\Psi}_{\lambda,x}) \subseteq WF_A(\Psi_{\lambda,x})\sqcup WF_A(\overline{\Psi}_{\lambda,x})$. Without loss of generality, suppose that $(y,\xi) \in WF_A(\Psi_{\lambda,x})$ but $(y,\xi) \notin WF_A(\Theta_{\pm})$. We can write $\Psi_{\lambda,x} = \Theta_{\pm} \mp \overline{\Psi}_{\lambda,x}$. It follows that $WF_A(\Psi_{\lambda,x}) \subseteq WF_A(\Theta_{\pm}) \cup WF_A(\overline{\Psi}_{\lambda,x})$. Since $(y,\xi) \notin  WF_A(\Theta_{\pm}) \cup WF_A(\overline{\Psi}_{\lambda,x})$ by assumption, this implies that $(y,\xi) \notin  WF_A(\Psi_{\lambda,x})$, which is a contradiction. Using the same argument, we obtain that if $(y,\xi) \in WF_A(\overline{\Psi}_{\lambda,x})$, then $(y,\xi) \in WF_A(\Theta_{\pm})$, which proves the corollary.
   
\end{proof}

Since the wavefront set of the distribution $\Psi_{\lambda,x}$ satisfies $WF_A(\Psi_{\lambda,x}) \cap -WF_A(\Psi_{\lambda,x}) = \emptyset$ for all $x \in \dS$ and the same holds for $\overline{\Psi}_{\lambda,x}$, we obtain the following corollary.

\begin{corollary}
The distributions $\Psi_{\lambda,x}$ and $\overline{\Psi}_{\lambda,x}$ cannot vanish on any open region of $\dS$.
\end{corollary}

Let $(L, M)$ be a symmetric space and $\Theta$ be a distribution on $L/M$. Then $L$ acts on $\Theta$ by 
\begin{align*}
	\pi_{-\infty}(g)\Theta (\varphi) = \Theta(\pi_{\infty}(g^{-1})\varphi), \qquad \varphi \in \mathcal{D}(L/M),
\end{align*}
where $\pi_{\infty}(g)\varphi(x) = \varphi(g^{-1} x)$.
\begin{definition}
	We say that a distribution $\Theta$ is $M$-invariant if $\pi_{-\infty}(h)\Theta = \Theta$ for all $h \in M$.
\end{definition}
\begin{definition}
	A distribution $\Theta$ is said to be a spherical distribution if it is an $M$-invariant  eigendistribution of the Laplace-Beltrami operator $\Delta$ on $L/M$. The space of all $M$-invariant spherical distributions satisfying  $\Delta \Theta = (\rho^2-\lambda^2)\Theta$ is denoted by $\mathcal{D}_\lambda^{*}(L/M)^M$.
\end{definition}

Consider $G' = \mr{O}_{1,n}(\R{})$, $H' = \mr{O}_{1,n-1}(\R{})$, $G = \mr{SO}_{1,n}(\R{})_e$ and $ H = \mr{SO}_{1,n-1}(\R{})_e$. Let $\mathcal{D}_{\lambda}^{*}(\dS)^{H'}$ and $\mathcal{D}_{\lambda}^{*}(\dS)^H$  be the spaces of $H'$ and $H$-invariant spherical distributions, respectively, on the de Sitter space satisfying $\Delta \Theta = (\rho^2- \lambda^2) \Theta$.

\begin{theorem} We have
\begin{enumerate}[leftmargin=*,label={\arabic*)}]
\item The dimension of $\mathcal{D}_{\lambda}^{*}(\dS)^H$ is 4.

\item The dimension of  $\mathcal{D}_{\lambda}^{*}(\dS)^{H'}$ is 2.
\end{enumerate}
\end{theorem}

Part (1) of the above theorem is \cite[Thm. 9.2.5]{D09} and {part (2) follows from \cite{OSe80}.}

We will now restate and prove Theorems~\ref{thm:2} and ~\ref{thm:3}. Let $\Psi_\lambda^\pm = \Psi_{\lambda,\pm e_n}$ and $\oPsi_\lambda^\pm = \oPsi_{\lambda,\pm e_n}$ as in the introduction.

\begin{manualtheoremx}{2}  Let $n \geq 2$ and let $\lambda\in \C{} \setminus \{(\rho + \N_0) \cup (-\rho-\N_0)\}$. Then the following holds:
	\begin{enumerate}[leftmargin=*,label={\arabic*)}]
	\item The distributions $\Psi_{\lambda}^+ $, $\Psi_{\lambda}^- $, $\oPsi_{\lambda}^+  $ and $\oPsi_{\lambda}^- $  are spherical 
	and form a basis for the space $\DdSl^H$.
	
	\item The distributions in part (1) 
	% $\Psi_{\lambda, e_n} $, ${\oPsi}_{\lambda, e_n} $, $\Psi_{\lambda, -e_n}$ and $\oPsi_{\lambda, -e_n}$
	are 
	positive definite if and only if $\lambda \in i\R{} \cup (-\rho,\rho)$. 
	
	\item Let $\Theta \in\DdSl^H$, $\Theta\not= 0$. Then
	\[ WF_A(\Theta)\subset WF_A(\Psi_{\lambda}^+) \sqcup WF_A( \oPsi_{\lambda}^+) \sqcup 
	WF_A(\Psi_{\lambda}^- )\sqcup WF_A (\oPsi_{\lambda}^-). \]
	Furthermore, if    $WF_A(\Theta)$ coincides with the analytic wavefront set of any of the distributions in part (1),  then $\Theta$ is
	a constant multiple of that distribution.
	\end{enumerate}  
\end{manualtheoremx}
\begin{proof}  We obtain from Lemma~\ref{lemma:psi} that all of them are eigendistributions of $\Delta$. From Theorem~\ref{thm:1} it follows that $\Psi_\lambda^+$, $\oPsi_\lambda^+$, $\Psi_\lambda^-$ and $\oPsi_\lambda^-$ are $H$-invariant and have disjoint analytic wavefront sets. Now, it follows from the definition of analytic wavefront sets that for some distributions $\Theta_1$ and $\Theta_2$ and for some constants $c_1$ and $c_2$,
	\begin{equation}\label{eq:wvadd}
		WF_A(c_1\Theta_1 + c_2\Theta_2) \subset WF_A(\Theta_1) \cup WF_A(\Theta_2).
	\end{equation}
   Let $\Theta = a\Psi_\lambda^+ + b\oPsi_\lambda^+ + c\Psi_\lambda^- + d\oPsi_\lambda^-$. Suppose that $\Theta = 0$. Then $WF_A(0) = \emptyset$. Without loss of generality, assume that $d \neq 0$. Then $\oPsi_\lambda^- = \frac{1}{d}(\Theta - a\Psi_\lambda^+ -  b\oPsi_\lambda^+ - c\Psi_\lambda^-)$. From (\ref{eq:wvadd}), we obtain that
	$WF_A(\oPsi_\lambda^-) \subset WF_A(\Theta) \cup WF_A(\Psi_\lambda^+) \cup  WF_A(\oPsi_\lambda^+) \cup WF_A (\Psi_\lambda^-)$.  Thus, we arrive at a contradiction, as $WF_A(\oPsi_\lambda^-)$ is disjoint from $WF_A(\Psi_\lambda^+) \cup  WF_A(\oPsi_\lambda^+) \cup WF_A (\Psi_\lambda^-)$. This implies that $d = 0$. Similarly, we obtain that $a= b= c = 0$. Therefore, these four distributions are linearly independent and span $\mathcal{D}_{\lambda}^*(\dS)^H$. This proves part (1).
	
Part (2) follows from  Theorem~\ref{thm:kernel}. When $\lambda \in i\R{} \cup (-\rho, \rho)$, the kernels $\Psi_\lambda$ and ${\oPsi}_{\lambda}$ are positive definite. Consequently, their limits are also positive definite.

Suppose that $\Theta \neq 0$ and  $WF_A(\Theta) = WF_A(\Psi_\lambda^+)$. Write $\Theta = a\Psi_\lambda^+ + b\oPsi_\lambda^+ + c\Psi_\lambda^- + d\Psi_\lambda^-$ and let $d \neq 0$. Part (3) follows from (\ref{eq:wvadd}). Following the same arguments as in part (1), we arrive at a contradiction and obtain $d = 0$. Similarly, we obtain $b = c = 0$. Finally, we have $a \neq 0$, which is a consequence of the assumption that $WF_A(\Theta) = WF_A(\Psi_\lambda^+) \neq \emptyset$.  By repeating the same arguments  when $WF_A(\Theta)$ coincides with analytic wavefront set of any of the other distributions, we establish the last claim. 
\end{proof}

\begin{manualtheoremx}{3}
	Let $n \geq 2$ and let $\lambda\in \C{} \setminus \{(\rho + \N_0) \cup (-\rho-\N_0)\}$. Then the following holds:
	\begin{enumerate}[leftmargin=*,label={\arabic*)}]
	\item The distributions $\Theta_{\lambda}^+ := \Psi_{\lambda}^+ + \oPsi_{\lambda}^+ $ and 
	$\Theta_{\lambda}^- := \Psi_{\lambda}^- + \oPsi_{\lambda}^-$  are $H^\prime$-spherical and form
	a basis for $\DdSl^{H^\prime} $.
	
	 \item The distributions in part  (1) are  positive definite if and only if $\lambda \in i\R{}\cup (-\rho,\rho)$.
	 
	 \item Let $\Theta \in \DdSl^{H^\prime}$, $\Theta \not=0$. Then
	 \[WF_A(\Theta)\subset WF_A(\Theta_{\lambda}^+)  \sqcup WF_A(\Theta_{\lambda}^-).\]
	 Furthermore, if
	 $WF_A(\Theta)$ coincides with the analytic wavefront set of any of the distributions in part (1), then  $\Theta$ is 
	 a constant multiple of that distribution. 
	\end{enumerate}
\end{manualtheoremx}
\begin{proof} Let 
	\[\Lambda_1 = \begin{pmatrix}
		-1 & 0 \\ 0 & I_{n}
	\end{pmatrix}, \qquad
	\Lambda_2 = \begin{pmatrix}
		1 & 0 & 0\\ 
		0 & -1 & 0 \\
		0 & 0 & I_{n-1}
	\end{pmatrix} .\]
	Then it is well known that \[\mr{O}_{1,n}(\R{}) = \mr{SO}_{1,n}(\R{})_e \sqcup  \Lambda_1 \mr{SO}_{1,n}(\R{})_e \sqcup \Lambda_2 \mr{SO}_{1,n}(\R{})_e \sqcup \Lambda_1\Lambda_2 \mr{SO}_{1,n}(\R{})_e.\]
	Moreover, when $\mr{O}_{1,n-1}$ and $\mr{SO}_{1,n-1}$ are viewed as a subgroups of $G'$ and $G$, respectively, stabilizing the base point $e_n$ of $\dS$, we obtain $\mr{O}_{1,n-1}(\R{}) =$
	\[ \mr{SO}_{1,n-1}(\R{})_e \sqcup  \Lambda_1 \mr{SO}_{1,n-1}(\R{})_e \sqcup \Lambda_2 \mr{SO}_{1,n-1}(\R{})_e \sqcup \Lambda_1\Lambda_2 \mr{SO}_{1,n-1}(\R{})_e.\]
	 We obtain \[\pi_{-\infty}(\Lambda_1)(\Psi_\lambda^\pm) = \oPsi_\lambda^\pm \quad \text{and} \quad \pi_{-\infty}(\Lambda_1)(\oPsi_\lambda^\pm) = \Psi_\lambda^\pm.\] Also, 
	$\pi_{-\infty}(\Lambda_1\Lambda_2)(\Psi_\lambda^\pm) = \oPsi_\lambda^\pm$ and $\pi_{-\infty}(\Lambda_1\Lambda_2)(\oPsi_\lambda^\pm) = \Psi_\lambda^\pm$. Therefore, they are not $H'$-invariant. However, their sums $\Theta_\lambda^+$ and $\Theta_\lambda^-$ are $H'$-invariant. The rest follows by the same argument as in the proof of Theorem~\ref{thm:2}, using Corollary~\ref{cor: pmwf}.

\end{proof}

%--------------------------------------------------------------------------------------------------------------------------------------------------------------------------------------------------------------------------------------------------------------------------------------------------------------------------------------------------
\section{Proof of Theorem 1}\label{sec:proof}
 Let $X = \dS$ and $x = ge_n$ for some $g \in \mr{SO}_{1,n}(\R{})_e$.  We will identify $\rT_{e_n}(\dS)$ with the Minkowski space $\R{1,n-1}$ equipped with the bilinear form $[v,v] = -v_0^2+v_1^2...+ v_{n-1}^2$. The cotangent space at $e_n$ will be identified with the tangent space at $e_n$ as $\R{1,n-1}$. Let $\mathbb{L}_{n-1}  = \{ v \in \rT_{e_n}(\dS): -v_0^2+v_1^2+...+ v_{n-1}^2=0\}$ be the light cone in $\rT_{e_n}(\dS)$ and denote by $\mathbb{L}_{n-1}^*$ the light cone in the cotangent space $\rT_{e_n}^*(\dS)$.

\begin{proof}[{\bf Proof of Theorem 1}] We will first prove part \textit{(2)} of the theorem in several steps. We begin by calculating the analytic wavefront set of the distribution $ \Psi_{\lambda, e_n}$. In steps~ (1)--(4) we determine the analytic wavefront set of the pullback of the distribution $\Psi_{\lambda,e_n}$ to the local coordinate chart in $\rT_{e_n}(\dS)$. In step (5) we deduce the analytic wavefront set of this distribution on de Sitter space using Theorem~\ref{cw}. We again use Theorem~\ref{cw} to derive the analytic wavefront set of $\Psi_{\lambda,x}$ for any $x \in \dS$ in (6) and obtain the same for $\oPsi_{\lambda,x}$ in (7).   
\medskip

(1) Let $g = \mr{Id}$ and write $(U,\kappa) = (U_{e_n},\mr{Exp}_{e_n})$, where $U_{e_n} \subset \rT_{e_n}(\dS)$ is the local chart  around $e_n$ with the coordinate map
\[y= \kappa(v) := \mathrm{Exp}_{e_n}(v) = C([v,v])e_n + S([v,v])v\]
as in Lemma~\ref{lemma: exp_x}. Let $V = V_{e_n} := \kappa (U)$. By Lemma~\ref{lemma:psi}, the analytic singular support of $\Psi_{\lambda,e_n}$ lies on the  light cone of $e_n$ given by \[\text{sing supp} _A(\Psi_{\lambda,e_n}) = \{y \in \dS\colon y_n=1\}= \kappa(\mathbb{L}_{n-1}) \subset V.\] Thus, it suffices to compute the wavefront set in $V$. To this end, we first compute $WF_A(\kappa^*\Psi_{\lambda,e_n})$, where $\kappa^*\Psi_{\lambda,e_n}$ is the pullback of $\Psi_{\lambda, e_n}$ to $U$ under $\kappa$. We see that 
\[WF_A(\kappa^*\Psi_{\lambda,e_n}) \subset \mathbb{L}_{n-1} \times \R{1,n-1} \subset U \times \R{1,n-1}.\]

Furthermore, it follows from the discussion in  \cref{sec:diffoper} and Lemma~\ref{lemma:psi} that this distribution is a solution to $P\Theta = 0$, where  \[P=  -\frac{\partial^2}{\partial v_0^2} + \sum_{i=1}^{n-1} \frac{\partial^2}{\partial v_i^2} - (\rho^2-\lambda^2).\] As a consequence of Theorem~\ref{thm:char}, we have
$$WF_A(\kappa^*\Psi_{\lambda,e_n}) \subset \mathrm{Char} (P),$$
where $$\mathrm{Char} (P) = \{(v,\xi) \in U \times (\R{1,n-1}\setminus \{0\}): P_n(v,\xi) =0\}.$$
Here, $P_n$ is the principal symbol of the differential operator $P$ given by
\[P_n = \xi_0^2 - \sum_{i=1}^{n-1} \xi_i^2.\] Then 
\[{\rm Char} (P) = \{(v,\xi) : v \in U, \xi \in \mathbb{L}_{n-1}^*\},\]
where $\mathbb{L}_{n-1}^*$ is the light cone in the  cotangent space.
\medskip

(2) Consider the map $f:U \rightarrow \R{}$  defined by 
    $$ f(v) := \ltfrac{1+C([v,v])}{2}.$$ The tangent map of $f$ is
$$df_v = \frac{-S([v,v])}{4} \begin{pmatrix}-2v_0, &2 \bv \end{pmatrix}, \quad v \in U.$$
If $\eta \in \R{}$ and $v \in U $ satisfy ${}^tdf_v(\eta) = 0$, then either $\eta = 0$, $v=0$ or $S([v,v])=0$. Suppose that $S([v,v]) =0$. Using the relation $S(z^2) = \sin{z}/z$, we obtain $[v,v] = m^2\pi^2$ for $m \in \mathbb{Z}\setminus \{0\}$. Such a $v$ does not belong to the set $U$. Therefore, $S([v,v]) \neq 0$ for any $v \in U$. 

Now consider the restriction of the map $f$ to the open sets $U_{\pm} := U_{\pm v_0 > 0}$ and denote it by $f_{\pm}$. Thus, the normal set (defined in Theorem~\ref{cw}) of the restriction is 
\begin{align*}
	N_{f_\pm} & := \{(f_{\pm}(v), \eta) \in \R{} \times \R{}: {}^tdf_v(\eta) = 0\} \\
	&=\{(f_\pm(v), 0): v \in U_{\pm}\}.
\end{align*}
 Note that if $v=0$, then $df_v \equiv 0$. Hence, we study the two cases $v \neq 0$ and $v=0$ separately.  
\medskip

(3) In this step we determine the singularities at $v \neq 0$. We will write
$v = (v_0, \bv)$ and consider the distributions ${}_2F_1(x+i0)$ and ${}_2F_1(x -i0)$ as in Theorem~\ref{thm:hgf} and Example~\ref{eg:hgf}. It follows from Example~\ref{eg:hgf} that $N_{f_+} \cap WF_A({}_2F_1(x - i0)) = \emptyset$ and $N_{f_-} \cap WF_A({}_2F_1(x + i0)) = \emptyset$.  As a result of Theorem~\ref{cw}, the pullback $f_{\pm}^*({}_2F_1(x \mp i0))$ is well defined on $U_{\pm}$. 

Let $dv = \kappa^* d\mu_{\dS}(y)$ denote the pullback to $U$ of the restriction of the $G$-invariant measure $d\mu_{\dS}$ to $V$.    For $v_0>0$ and $\varphi \in C_c^\infty(U_{+})$, it follows that
\begin{align*}
    \kappa^*\Psi_{\lambda,e_n} (\varphi)	&= \lim_{t \nearrow \pi/2}\int_{U_+}\Psi_{\lambda}(z_t,\kappa(v)) \overline{\varphi}(v) dv \\
	&=  \lim_{t \nearrow \pi/2}\int_{U_+} {}_2F_1 \left( \frac{1+ [z_t, \kappa (v)]}{2}\right)  \overline{\varphi}(v) dv  \\
	&= {}_2F_1 \left( \frac{1+ C([v,v])}{2} - i0\right)({\varphi}) = {}_2F_1 (f_+(v)-i0) (\varphi)\\
	&= f_+^*({}_2F_1(x-i0)) (\varphi).
\end{align*}

By the same argument, we obtain $\kappa^*\Psi_{\lambda,e_n}(\varphi) = f_-^*({}_2F_1(x+i0)) (\varphi)$ for $\varphi \in C_c^\infty(U_-)$. From this we conclude that $\kappa^*\Psi_{\lambda,e_n} = f_{\pm}^*({}_2F_1(x\mp i0))$ on $U_+$ and $U_-$, respectively. As a consequence of Theorem~\ref{cw}, we obtain on $U_+$
\begin{equation*}
    \begin{aligned}
  WF_A(\kappa^*\Psi_{\lambda, e_n}) &=  WF_A({}_2F_1(f_+(v) - i0)) \\ 
  & \subseteq \{(v, {}^tdf_{+,v}(\xi)): v_0 >0, (f_+(v), \xi) \in WF_A({}_2F_1(x - i0))\}\\
  & = \{ (v, \xi) : [v,v]=0,  \xi = (-v_0,\bv ), v_0 >0\},
\end{aligned}
\end{equation*}

    and on $U_-$
    \begin{equation*}
        \begin{split}
     WF_A(\kappa^*\Psi_{\lambda, e_n}) &= WF_A({}_2F_1(f_-(v) + i0))\\
     &\subseteq \{(v, {}^tdf_{-,v}(\xi)): v_0<0, (f_-(v), \xi) \in WF_A({}_2F_1(x + i0))\}\\
     &= \{ (v, \xi) :[v,v] =0, \xi = (v_0,-\textbf{v}), v_0 <0\}.
    \end{split}
    \end{equation*}
    
    From step (1), we know that the projection of $WF_A(\kappa^*\Psi_{\lambda, e_n})$ onto $U$ is \[\mr{proj}_U (WF_A(\kappa^*\Psi_{\lambda, e_n}) ) = \text{sing supp}_A (\kappa^*\Psi_{\lambda, e_n}) = \{v \in U\: [v,v]=0\}.\]
     Hence, for every $v$ in the above set that lies in $U_{\pm}$, the fiber of $WF_A(\kappa^*\Psi_{\lambda, e_n})$ at $v$ must be non-empty. Since the analytic wavefront set is conic and the inclusions obtained above  allow only one direction at $v$, the inclusions must be equalities. More precisely, in the domain $U_+$ we have
     \begin{equation}\label{eq:a}
     	WF_A(\kappa^*\Psi_{\lambda, e_n})  = \{ (v, \xi) : [v,v]=0,  \xi = (-v_0,\bv ), v_0 >0\},
     \end{equation}
   
     and in $U_-$ we have
     \begin{equation}\label{eq:b}
     	WF_A(\kappa^*\Psi_{\lambda, e_n})  =\{ (v, \xi) :[v,v] =0, \xi = (v_0,-\textbf{v}), v_0 <0\}.
     \end{equation}

\medskip

(4) Next, we determine the singularities at $v=0$. For this we will use the propagation of singularities theorem (see Theorem~\ref{PS}), which states that the analytic wavefront set is invariant under the Hamiltonian flow of $P_n$ whenever ${\partial P_n (v,\xi)}/{\partial \xi} \neq 0$. 
We have ${\partial P_n}(v, \xi)/{\partial \xi} = 2(\xi_0, -\xi_1,..., -\xi_{n-1}) = 0$  only if $\bf{\xi} =0$. As $\xi = {\bf 0}$ is excluded from the definition of analytic wavefront sets, we can apply \cref{PS} at $v=0$.
We know that \[\begin{split}
	\mr{Char} (P_n) &= \left\{(v,\xi) \in \rT^*(U)\: \xi_0^2-\sum_{i=1}^{n-1}\xi_i^2  =0  \right\} \\
	&= \{(v,\xi): v\in U, \xi \in L_{n-1}^*\}.
\end{split}\]
 Let $(v(t),\xi(t))$ be the Hamiltonian strip of $P_n$ with initial data $v(0)=v^0$ and $\xi(0) = \xi^0 \in \mathbb{L}_{n-1}^*$. Then we must have $\xi(t) \in \mathbb{L}_{n-1}^*$ for all $t \in \R{}$ and 
 $$\frac{\partial v}{\partial t} = \frac{\partial P_n}{\partial \xi}, \qquad \frac{\partial \xi}{\partial t} = -\frac{\partial P_n}{\partial v}.$$
That is, the following relations hold:
\begin{equation*}\label{eq:hamilton}
    \begin{aligned}
    \dot{v_0} = 2\xi_0 &, \qquad \dot{\xi_0}=0 \\
    \dot{v_i} = -2\xi_i&, \qquad \dot{\xi_i} = 0 \quad \text{for}\; i=1,...,n-1.
\end{aligned}
\end{equation*}
This gives
\begin{equation}\label{eq:Ham1}
\xi(t) = \xi^0,
\end{equation}
which in turn gives
\begin{equation}\label{eq:Ham2}
\begin{split}
v_0(t) &=v^0_0+ 2\xi_0^0t, \\
v_i(t) &= v^0_i -2\xi_i^0t \quad \text{for}\; i=1,...,n-1.
\end{split}	  
\end{equation}

We now show that $\{(0,\xi): \xi \in \mathbb{L}_{n-1}^*, \xi_0 < 0\} \subseteq WF_A(\kappa^*\Psi_{\lambda,e_n})$. Let $\xi^0 \in \mathbb{L}_{n-1}^*$ with $\xi_0 < 0$. Let $v^0 \in U$ such that $\xi_0^0 = -v_0^0 < 0$, $\xi_i^0 = v_i^0$. From (\ref{eq:a}) we have $(v^0,\xi^0) \in WF_A(\kappa^*\Psi_{\lambda,e_n})$ with $v^0 \in U_+$. Then the Hamiltonian strip with the above initial data $(v^0,\xi^0)$ is given by $v_0(t) = v_0^0(1-2t)$, $v_i(t) = v_i^0(1-2t)$,  $\xi_0(t) = -v_0^0 < 0$, and $\xi_i(t) = v_i^0$, for all $i=1,\ldots, n-1$. The curve $v(t)$ passes through the origin at $t = 1/2$.  As a result of \cref{PS}, $(v(1/2),\xi(1/2)) =(0,(-{v_0}^0, {\textbf{v}}^0))$ must be in the wavefront set, establishing the above inclusion. We obtain the same using (\ref{eq:b}).

On the other hand, we now claim that the fiber over $0$ in the wavefront set of $\kappa^*\Psi_{\lambda, e_n}$ is precisely the set $\{(0,\xi): \xi \in \mathbb{L}_{n-1}^*, \xi_0 < 0\}$. Suppose $(0,\xi^0) \in WF_A(\kappa^*\Psi_{\lambda, e_n})$ with $\xi_0^0 \geq 0$. Then, using (\ref{eq:Ham1}) and (\ref{eq:Ham2}), the Hamiltonian strip with initial data $(0,\xi^0)$ is given by $v_0(t) = 2\xi_0^0t$, $v_i(t) = -2\xi_i^0t$, $\xi_0(t) = \xi^0_0$, and $\xi_i(t) = \xi^0_i$, for $i=1,\ldots,n-1$. 

Again by \cref{PS}, we must have $(v(t),\xi(t))$ in the wavefront set of $\kappa^*\Psi_{\lambda, e_n}$ for all $t \in \R{}$.
Without loss of generality, let $t=1$.
Observe that $[v(1),v(1)] = 4[\xi^0,\xi^0] = 0$ with $v_0(1) > 0$. Then $v(1) \in U_+$ and $(v(1),\xi(1))$ lies in the left hand side of (\ref{eq:a}). This must imply that $\xi_0(1)<0$, contradicting the assumption that $\xi_0(1) = \xi_0^0 >0$. Therefore, we have established the claim. 

More precisely, by combining steps (3) and (4), we obtain
\begin{small}
    \begin{equation}\label{eq:WFkappa}
    \begin{aligned}
WF_A(\kappa^*\Psi_{\lambda, e_n}) 
  &= \{(0,\xi): \xi_0 <0, \xi \in \mathbb{L}_{n-1}^* \} \cup \{(v, \xi): \xi = (-v_0,\textbf{v}), v_0>0\} \\
  &\cup \{(v,\xi): \xi = (v_0,-\textbf{v}), v_0 < 0 \},
\end{aligned}
\end{equation}
\end{small}
for $v \in \mathbb{L}_{n-1}^*$.
\medskip

(5) In this step we determine the wavefront set of $\Psi_{\lambda,e_n}$ on de Sitter space using the pullback of the set $WF_A(\kappa^*\Psi_{\lambda,e_n})$ to $V$ under the map $(\kappa^{-1})^*$. If $v \in \mathbb{L}_{n-1}$, then $y= C[v,v]e_n + S[v,v]v = e_n + v$ and $[y-e_n,y-e_n]=0$. This implies that $y$ lies on the light cone of $e_n$ on $\dS$. Using Theorem~\ref{cw} and (\ref{eq:WFkappa}), we conclude that the analytic wavefront set of $\Psi_{\lambda, e_n}$ is given by 
$$WF_A(\Psi_{\lambda, e_n}) = (\kappa^{-1})^*WF_A(\kappa^*\Psi_{\lambda, e_n}).$$
By identifying $\mathbb{L}_{n-1}^*$ with $\mathbb{L}_{n-1}$, we get, for all $v \in \mathbb{L}_{n-1}$ 
\begin{align*}
  WF_A(\Psi_{\lambda, e_n}) &=  \{(e_n,v): v_0>0 \} \cup \{(e_n +v,  (-v_0, {\bf v})): , v_0>0\}\cup \\ 
  &\{(e_n+v, (v_0, -{\bf v})): v_0 >0\}.
\end{align*}
\medskip

(6) Consider $\Psi_{\lambda,x} = \lim_{t \nearrow \pi/2}\Psi_\lambda(g z_t,y)$. For $x = g  e_n$, let $\ell_g : \dS \rightarrow \dS$ be the left translation given by $\ell_g(y) = g y$. Since $\ell_g$ is an analytic diffeomorphism, $d\ell_g$ is an isomorphism of tangent spaces. Therefore, $N_{\ell_g} = \{(y, 0): y \in \dS\}$. Hence, by Theorem~\ref{cw}, the pullback of the distribution $\Psi_{\lambda, e_n}$ under the map $\ell_g^{-1}$ is $\Psi_{\lambda,x}$ and  
 $WF_A(\Psi_{\lambda,x}) = (\ell_g^{-1})^* WF_A(\Psi_{\lambda, e_n})$. That is, $(y,\xi) \in WF_A(\Psi_{\lambda,x})$ if and only if $(\ell_g^{-1}(y), {}^td\ell_{g}^{-1}\xi) \in WF_A(\Psi_{\lambda, e_n})$. This implies that $y = x + gv$ for some $v \in \mathbb{L}_{n-1}$. As $G$ acts transitively on the light cone, preserving the sign of the first coordinate, ${}^td\ell_{g^{-1}}\xi =  g^{-1} \xi$ also lies in $\mathbb{L}_{n-1}^*$ with $\sgn{({}^td\ell_{g^{-1}}\xi)_0} = \sgn{\xi_0}$.
 
It follows that
\begin{align*}
WF_A(\Psi_{\lambda,x}) &= \{(x,v): v_0 < 0 \} \cup \{(x +v,  (-v_0, {\bf v})): v_0<0\}\cup \\ 
  &\{(x+v,(v_0, -{\bf v})): v_0 <0\},
\end{align*}
for $v \in \rT_x(\dS)$ such that $[v,v]=0$
\medskip

(7) Analogously, we obtain   $\kappa^*\overline{\Psi}_{\lambda, e_n}(\varphi) = f_+^*({}_2F_1(x + i0)) (\varphi)$ for $ \varphi \in C_c^\infty(U_+)$ and $\kappa^*\overline{\Psi}_{\lambda, e_n}(\varphi) = f_-^*({}_2F_1(x - i0)) (\varphi)$ for $\varphi \in C_c^\infty( U_-)$. Following the same steps as above, we obtain
\begin{align*}
  WF_A(\overline{\Psi}_{\lambda, e_n}) &=  \{(e_n,v): v_0 > 0 \} \cup \{(e_n +v,  (v_0, -{\bf v})): v_0>0\} \\ 
  &\cup\{(e_n+v, (-v_0, {\bf v})): v_0 >0\},
\end{align*}
for $v \in \rT_{e_n}(\dS)$ with $[v,v]=0$ and
\begin{align*}
   WF_A(\overline{\Psi}_{\lambda,x}) &= \{(x,v): v_0 > 0 \} \cup \{(x +v, (v_0, -{\bf v})): v_0>0\} \\ 
  &\cup\{(x+v, (-v_0, {\bf v})): v_0 >0\},
 \end{align*}
 for $v \in \rT_{x}(\dS)$ with $[v,v]=0$.
Hence, we conclude part \textit{(2)} of the theorem.

For part \textit{(1)}, let $x = g_1e_n=g_2e_n$. We denote $\Theta_1 = \lim_{t \nearrow \pi/2}\Psi_{\lambda}(g_1z_t,\cdot)$  and $\Theta_2 = \lim_{t \nearrow \pi/2}\Psi_{\lambda}(g_2z_t,\cdot)$. From the proof of Lemma~\ref{lemma:psi}, we have $\Theta_1(y) = \Theta_2(y)$ in the open regions $\Gamma^+(x), \Gamma^-(x)$ and $(\overline{\Gamma(x)})^c$. Moreover, observe that $(\ell_{g_i}^{-1})^* \Psi_{\lambda,e_n} = \Theta_i$, for $i=1,2$. From Theorem~\ref{cw}, we have \[WF_A(\Theta_i) = (\ell_{g_i}^{-1})^* WF_A(\Psi_{\lambda,e_n}), \qquad i=1,2.\]
Note that for $v \in \rT_{e_n}(\dS)$ with $[v,v]=0$, we have $g_iv \in \rT_x(\dS)$ with $[g_iv,g_iv] = 0$. Thus, we obtain
$WF_A(\Theta_1) = WF_A(\Theta_2) = WF_A(\Psi_{\lambda,x})$. Then an application of Theorem~\ref{thm:supp} to the distribution $\Theta_1-\Theta_2$ gives $\Theta_1 = \Theta_2$, proving part \textit{(1)} of the theorem.
\end{proof}

%--------------------------------------------------------------------------------------------------------------------------------------------------------------------------------------------------------------------------------------------------------------------------------------------------------------------------------------------------
\appendix 

 \section{Boundary values of the hypergeometric function}\label{sec:hgf}

We write ${}_2F_1(a,b;c,z) = {}_2F_1(z)$ as before. In this section, we show that ${}_2F_1(x+i0) := \lim_{y \searrow 0}{}_2F_1(x+iy)$ and ${}_2F_1(x-i0) := \lim_{y \searrow 0}{}_2F_1(x-iy)$ are distributions for $a = \rho + \lambda$, $b = \rho - \lambda$ and $c= n/2$. It is a fact that $\hgf{z}$ has a branch cut on $[1,\infty)$. Note that the convergence for $x < 1$ is uniform on compact sets. Given two
nonnegative functions $f$ and $g$ on a domain D, we write $f \asymp g$ if there exist positive
constants $C_1$ and $C_2$ such that $C_1g(x) \leq f(x) \leq C_2g(x)$ for all $x \in D$.

\begin{lemma}\label{lemma:hgf}
	Let $\lambda \in \C{} \setminus \{(\rho + \mathbb{N}_0)\cup (-\rho - \mathbb{N}_0)\}$. The limits \[ \lim_{y \searrow 0}\; {}_2F_1(\rho + \lambda, \rho - \lambda, n/2, x + iy), \quad \lim_{y \searrow 0}\; {}_2F_1(\rho + \lambda, \rho - \lambda, n/2, x - iy) \] converge uniformly on compact subsets of the regions $x>1$ and $x<1$, where they are analytic. On the interval $1<x<2$, the limits are given by the real analytic functions as follows:  
	\begin{equation}\label{eq:2}
		\begin{split}
			&\hgf{x\pm i0} =  \frac{\Gamma(n/2)\Gamma((2-n)/2)}{\Gamma(1/2 + \lambda)\Gamma(1/2-\lambda)} {}_2F_1(\rho + \lambda,\rho-\lambda; \frac{n}{2};1-x) \\
			&+ e^{\mp i\pi(\frac{2-n}{2})}(x-1)^{\frac{2-n}{2}}  \ltfrac{\Gamma(n/2)\Gamma((n-2)/2)}{\Gamma(\rho + \lambda)\Gamma(\rho-\lambda)} {}_2F_1(1/2- \lambda,1/2+\lambda; \frac{4-n}{2};1-x),
		\end{split}
	\end{equation} 
	if $n$ is odd and  $\hgf{x \pm i0} =$

		\begin{equation}\label{eq:2'} 	
			\begin{split}
				(-1)^{\frac{2-n}{2}}& \Bigg( \ltfrac{\Gamma(n/2)}{\Gamma(\rho + \lambda)\Gamma(\rho-\lambda)}\sum_{k=0}^{n/2-2}\ltfrac{(n/2-k-2)!(1/2+\lambda)_k(1/2-\lambda)_k}{k!}(x-1)^{k+1-\frac{n}{2} }\\
				&+ \ltfrac{\Gamma(n/2)}{\Gamma(1/2 + \lambda)\Gamma(1/2-\lambda)}\sum_{k=0}^{\infty}\ltfrac{(\rho+\lambda)_k(\rho-\lambda)_k}{k!(n/2 - 1 + k)!}[\psi(k+1) + \psi(n/2 + k) \\
				& - \psi(\rho+\lambda +k) - \psi(\rho-\lambda +k) - \ln (x-1) \pm i\pi](-1)^k(x-1)^{k} \Bigg),
			\end{split} 
		\end{equation}
		if $n$ is even. Here, $\psi(z) = \Gamma'(z)/\Gamma(z)$. Furthermore, the asymptotic behavior near $z=1$ is given by  
	\begin{equation}\label{eq:3}
		{}_2F_1(z) \asymp  \ltfrac{1}{\Gamma(\rho+\lambda)\Gamma(\rho-\lambda)} (-\ln{(1-z}))
	\end{equation}
	for $n = 2$, and 
	\begin{equation}\label{4}
		\hgf{z} \asymp \frac{\Gamma(n/2)\Gamma((n-2)/2)}{\Gamma(\rho+\lambda)\Gamma(\rho-\lambda)}(1-z)^{\frac{2-n}{2}}
	\end{equation}
	for $n \geq 3$.
\end{lemma}
\begin{proof}
	Let $n \geq 2$. If $n$ is odd, then 
	$c-a-b = \frac{2-n}{2}$ is not an integer.
	Therefore, for $|z-1|<1$ and $|\mr{arg}(1-z)|<\pi$, we can use the following transformation 
	\begin{equation*}\label{eq: 1}
		\begin{split}
			&\hgf{z} = \ltfrac{\Gamma(n/2)\Gamma((2-n)/2)}{\Gamma(1/2 + \lambda)\Gamma(1/2-\lambda)} {}_2F_1(\rho + \lambda,\rho-\lambda; \frac{n}{2};1-z) \\
			&+ (1-z)^{\frac{2-n}{2}}  \ltfrac{\Gamma(n/2)\Gamma((n-2)/2)}{\Gamma(\rho + \lambda)\Gamma(\rho-\lambda)} {}_2F_1(1/2- \lambda,1/2+\lambda; \frac{4-n}{2};1-z).
		\end{split}
	\end{equation*}
	Suppose that $1<x<2$. Using the above transformation, we obtain the expression (\ref{eq:2}).
	For $x\geq 2$, we can use linear transformations of hypergeometric functions to extend $\hgf{x\pm i0}$ analytically. If $n$ is even, we obtain expression (\ref{eq:2'}) for $1 < x <2$ from \cite[ (9.7.5), (9.7.6)]{LS}. Similarly, we can extend $\hgf{x \pm i0}$ analytically for $x \geq 2$ using the formulae from \cite[Sec. 9.7]{LS}.
	
	To deduce the asymptotic behavior near $z=1$, we use Euler's transformation \cite[(9.5.3)]{LS} for $|\arg (1-z)|<\pi$, which gives
	\[\hgf{\rho+\lambda,\rho-\lambda;n/2;z} =(1-z)^{\frac{2-n}{2}}\hgf{1/2-\lambda,1/2+\lambda;n/2;z}.\]  
	Using the Taylor series expansion of $\hgf{1/2-\lambda,1/2+\lambda;n/2;z}$ around $z=1$ with $|\arg (1-z)|<\pi$, we obtain
	\[\hgf{\rho+\lambda,\rho-\lambda;n/2;z} \asymp \frac{\Gamma(n/2)\Gamma((n-2)/2)}{\Gamma(\rho+\lambda)\Gamma(\rho-\lambda)}(1-z)^{\frac{2-n}{2}}\]
	for $n\geq 3$. For $n=2$, using \cite[(9.7.5)]{LS}, we obtain
	\[\hgf{1/2+\lambda,1/2-\lambda;1;z}\asymp  \ltfrac{1}{\Gamma(\rho+\lambda)\Gamma(\rho-\lambda)} (-\ln{(1-z})).\]
	Moreover,
	from (\ref{eq:2}), (\ref{eq:2'}) and \cite[(9.3.4)]{LS} we obtain for $n\geq 3$
	
	\begin{equation*}
		\underset{x \rightarrow 1^+}{\mathrm{lim}} \frac{\hgf{x\pm i0}}{(x-1)^{c-a-b}} =  e^{\mp i \pi (\frac{2-n}{2})}\frac{\Gamma(n/2)\Gamma((n-2)/2)}{\Gamma(\rho +\lambda)\Gamma(\rho-\lambda)}.
	\end{equation*}
	For $n=2$, the first summation in (\ref{eq:2'}) does not appear. Thus, we obtain
	
	\begin{equation*}
		\underset{x \rightarrow 1^+}{\mathrm{lim}} \frac{\hgf{x \pm i0}}{- \ln (x-1) \pm i \pi} = \frac{1}{\Gamma(\rho + \lambda)\Gamma(\rho-\lambda)}.\label{eq:6}
	\end{equation*}
	This proves the lemma.
\end{proof}

\begin{theorem} \label{thm:hgf}
	Let $\lambda \in \C{} \setminus \{(\rho + \mathbb{N}_0)\cup (-\rho - \mathbb{N}_0)\}$. For $y>0$, the limits \[ \lim_{y \searrow 0}\; {}_2F_1(\rho + \lambda, \rho - \lambda, n/2, x + iy), \quad \lim_{y \searrow 0}\; {}_2F_1(\rho + \lambda, \rho - \lambda, n/2, x - iy) \] exist in the sense of distributions. Moreover, the limits converge uniformly on compact subsets of the regions $x>1$ and $x<1$, where the distributions are represented by analytic functions.

\end{theorem}
\begin{proof}
Let $z=x+iy$. We know from Lemma~\ref{lemma:hgf} that the limits converge uniformly on compact sets in the regions $x>1$ and $x<1$.
Moreover, the growth of $\hgf{z}$ near $z = 1$ for $n \geq 3$ is $$|\hgf{x \pm iy}| \leq C_1 \Big|\ltfrac{\Gamma(n/2)\Gamma((n-2)/2)}{\Gamma(\rho + \lambda)\Gamma(\rho - \lambda)}\Big|(|1-z|^{\frac{2-n}{2}}) \leq C_2 |y|^{\frac{2-n}{2}},$$
where  $C_1,C_2 >0$ are constants. It follows from \cite[Thm. 3.1.11]{H63} that the limits converge to distributions.  

Similarly, for $n=2$, there exists a constant $C>0$ such that
\[ \hgf{x \pm iy} \leq C \left|\frac{1}{\Gamma(\rho + \lambda)\Gamma(\rho-\lambda)} ( \ln (1-z))\right|\]
for $|1-z|<\epsilon$, where $\epsilon>0$ is small enough. Since the logarithm is integrable on compact sets, it follows that $\hgf{x\pm i0}$ is a distribution.
\end{proof}

%----------------------------------------------------------------------------------------------------------------------------------------------------------------------------------------------------------------------------------------------------

\section{{Distributions: $(x\pm i0)^\frac{2-n}{2}$ and $\ln (x\pm i0)$}} \label{appendixa}
We recall the distributions $(x\pm i0)^{\frac{2-n}{2}}$ and $\ln (x \pm i0)$ (see \cite{GS64}). 
Let $n$ be odd and consider the distributions $x^{\frac{2-n}{2}}_+$ and $x^{\frac{2-n}{2}}_-$. Let $\varphi \in C_c^\infty(\R{})$. For $n=3$,
\begin{equation*}
    (x^{-\frac{1}{2}}_+, \varphi) = \int_0^\infty x^{-\frac{1}{2}}\varphi(x) dx
\end{equation*}

defines a regular distribution. However, for $n \geq 5$ and $m= (n-5)/2$, we have

\begin{equation*}
    (x^{\frac{2-n}{2}}_+, \varphi) = \int_0^\infty x^{\frac{2-n}{2}}\Big[\varphi(x) - \varphi(0)-x\varphi'(0)-...-\ltfrac{x^m}{(m)!}\varphi^m(0) \Big] dx.
\end{equation*}

The distribution $x^{\frac{2-n}{2}}_-$ is defined as follows:
\begin{equation*}
    (x^{\frac{2-n}{2}}_-, \varphi(x)) = (x^{\frac{2-n}{2}}_+, \varphi(-x)).
\end{equation*}

Let $n$ be even. Set $k = (n-2)/2$. If $k$ is even,

\begin{equation*}
\begin{aligned}
    (x^{-k}, \varphi) &= \int_0^\infty x^{-k} \Big( \varphi(x) + \varphi(-x) \\
    -2\Big[\varphi(0) &+ \frac{x^2}{2!}\varphi''(0)+...+\frac{x^{k-2}}{(k-2)!}\varphi^{k-2}(0) \Big]\Big) dx.
    \end{aligned}
\end{equation*}

If $k$ is odd,

\begin{equation*}\label{eq:nEvenkOdd}
\begin{aligned}
    (x^{-k}, \varphi) &= \int_0^\infty x^{-k} \Big( \varphi(x) - \varphi(-x) \\
    -2\Big[x\varphi'(0) &+ \frac{x^3}{3!}\varphi'''(0)+...+\frac{x^{k-2}}{(k-2)!}\varphi^{k-2}(0) \Big]\Big) dx.
    \end{aligned}
\end{equation*}

Consider the following distributions:
\begin{equation*}
    (x \pm i0)^{\frac{2-n}{2}} = \lim_{y \searrow 0} (x \pm iy)^{\frac{2-n}{2}}.
\end{equation*}

When $n$ is odd, 
\begin{align*}
     (x + i0)^{\frac{2-n}{2}} &=  x^{\frac{2-n}{2}}_+ + e^{i \pi \frac{2-n}{2}} x^{\frac{2-n}{2}}_-,\\
     (x - i0)^{\frac{2-n}{2}} &=  x^{\frac{2-n}{2}}_+ + e^{-i \pi \frac{2-n}{2}} x^{\frac{2-n}{2}}_-.
\end{align*}
   
 Let $n$ be even and set $k = (n-2)/2$. Then we have 
\begin{align*}
     (x + i0)^{-k} &=  x^{-k} - \ltfrac{i \pi (-1)^{k-1}}{(k-1)!}\delta^{k-1}(x),\\
    (x - i0)^{-k} &=  x^{-k} + \ltfrac{i \pi (-1)^{k-1}}{(k-1)!}\delta^{k-1}(x).
\end{align*}

Finally, we have the distributions 
\begin{equation*}
    \ln (x \pm i0) = \lim_{y \searrow 0} \ln (x \pm iy),
\end{equation*}
where

\begin{equation*}
    \ln (x \pm i0) = \begin{cases} \ln |x| \pm i\pi & \mathrm{for} \; x< 0,\\
    \ln x & \mathrm{for}\; x > 0.
    \end{cases}
\end{equation*}
%-----------------------------------------------------------------------------------------------------------------------------------------------------------------------------------------------------------------------------------------------------

\end{document}